\documentclass[11pt]{article}
\input epsf
\usepackage{amsmath}
\usepackage{amssymb}
\usepackage{amscd}
\newcommand{\qed}{\hbox{\unskip\nobreak\hfil
        \penalty50\hskip1em\hbox{}\nobreak\hfil
        $\square$\parfillskip=0pt\finalhyphendemerits=0 \par}}

\newtheorem{dfn}{Definition}[section]
\newtheorem{rem}[dfn]{Remark}
\newtheorem{Not}[dfn]{Notation}
\newtheorem{conv}[dfn]{Convention}
\newtheorem{thm}[dfn]{Theorem}
\newtheorem{lem}[dfn]{Lemma}
\newtheorem{prop}[dfn]{Proposition}
\newtheorem{cor}[dfn]{Corollary}

\newtheorem{defn}[dfn]{Definition}

\def\no{\noindent}

\def\={=\joinrel =}
\def\i{\sqrt[r]{I}}
\def\0{\emptyset}

\def\proof{\par\medskip\noindent{\it Proof: }}
\def\lra{\longrightarrow}

\def\>{\rangle}
\def\<{\langle}

\def\C{\mathbb C}
\def\t{\tilde}
\def\R{\mathbb R}
\def\P{\mathbb P}
\def\Z{\mathbb Z}
\def\E{{\cal E}}
\def\p{{\cal P}}

\def\a{{\cal A}}
\def\q{{\cal Q}}

\def\T{{\cal T}}
\def\V{{\cal V}}
\def\I{{\cal I}}
\def\A{{\mathbb A}}
\def\M{{\cal M}}
\def\eps{\epsilon}

\def\al{\alpha}
\def\be{\beta}

\def\Ga{\Gamma}

\def\Del{\Delta}
\def\Si{\Sigma}
\def\si{\sigma}
\def\L{{\cal L}}

\def\la{\lambda}

\def\8{\infty}

\def\M{{\cal M}}

\def\a{{\cal A}}
\def\b{{\cal B}}

\def\E{{\cal E}}
\def\O{{\cal O}}
\def\J{{\cal J}}

\def\mf{\mathfrak }

\def\nbd{neighborhood~}

\def\k{\hbox{\bf k}}

\def\1{\hbox{\bf 1}}

\def\u{\sqcup}
\def\v{\vee}

\def\BE{\begin{equation}}
\def\EE{\end{equation}}

\def\p{{\cal P}}
\begin{document}
\title{Universality theorems for configuration 
spaces of planar linkages}
\author{Michael Kapovich \ 
and John J.\ Millson}
\date{March 23, 1998}
\maketitle

\begin{abstract}
\no We prove realizability theorems for vector-valued polynomial mappings, 
real-algebraic sets and compact smooth manifolds by moduli spaces of planar 
linkages.  We also establish a relation between universality theorems for  
moduli spaces of mechanical linkages and projective arrangements.  
\end{abstract}

\section{Introduction} 

This paper deals with moduli spaces of planar linkages.  
An {\em abstract linkage} $(L, \ell)$ is a 
graph $L$ with a positive real number $\ell(e)$ assigned to 
each edge $e$. We assume that we have chosen a distinguished oriented edge
 $e^*= [v_1 v_2]$ in $L$. The {\em moduli space} $\M(\L)$ of planar 
realizations of $\L:= (L, \ell, e^*)$ is the set\footnote{Apriori this set 
could be empty.} of maps $\phi$ from the vertex set of $L$ into the 
Euclidean plane $\R^2$ (which will be identified with the complex plane 
$\C$) such that

\begin{itemize}
\item $|\phi(v)- \phi(w)|^2 = (\ell[vw])^2$ for each edge $[vw]$ of $L$. 
\item $\phi(v_1)=(0,0)$.
\item $\phi(v_2)= (\ell(e^*), 0)$. 
\end{itemize}

Clearly these conditions give $\M(\L)$ natural structure of a  
real-algebraic set in $\R^{2r}$ where $r$ is the number of vertices in $L$. 

It is important to note that the ideal $I$ in the polynomial ring 
$\R[X_1,Y_1,..., X_r,Y_r]$ (here $X_i, Y_i$ are the coordinates of 
$\phi(v_i)$) generated by the above equations can be strictly contained in 
the ideal of {\em all} polynomials vanishing on the set $\M(\L)\subset \R^{2r}$. 
We let ${\mf M}(\L)$ denote the affine subscheme of $\R^{2r}$ corresponding 
to the ideal $I$. We will in fact add more functions to the ideal $I$ 
corresponding 
to certain degenerate triangles -- see Convention \ref{contri}.

In Definition \ref{2.7} we define {\em functional linkages}. The abstract 
linkage corresponding to a functional linkage comes with two   
sets of vertices: 
the {\em inputs} $(P_1,..., P_m)$ and the {\em outputs} $(Q_1,..., Q_n)$. 
We let  $p: \M(\L) \to \A^m$ (the {\em input map}) and $q: \M(\L) \to \A^n$ 
(the {\em output map}) be the {\em forgetful  maps} that record only the 
positions of the images of $P_i$'s and $Q_j$'s  under realizations $\phi$. 
Here the affine line $\A$ is either $\C\cong \R^2$  (in which case we refer 
to $\L$ as a {\em complex functional linkage}) or 
$\R= \R \times  \{0\}\subset \R^2$ (in which case we refer 
to $\L$ as {\em real functional linkage}). 

We say that $\L$ as above is {\em functional}\, \footnote{See Definition 
\ref{2.7} for more precise definition.} for a mapping $f: \A^m \to \A^n$ 
if there is a commutative diagram

\[
\begin{array}{ccc}
~ & \M(\L) & ~ \\
~ & \swarrow p \quad  \quad q\searrow & ~  \\
\A^m & \stackrel{f}{\longrightarrow}  & \A^n\\
\end{array}
\] 

\noindent and $p$ is a {\em regular topological branched cover} of a 
bounded domain in $\A^m$. We prove

\bigskip
\noindent {\bf Theorem A.} {\em Let $f: \A^m \to \A^n$ be a 
polynomial map with real coefficients (where $\A$ is either $\C$ or $\R$). 
Let $\O$ be a point in  $\A^m$ and $r>0$. Then there is a functional linkage 
$\L$ for $f$ such that the $r$-ball $B_r(\O)$ is in the interior of the 
image of $p$ and $p$ is an analytically trivial covering over 
$B_r(\O)$. The same 
conclusion holds for polynomial maps whose coefficients are not required to be 
real if we use more general definition of functional linkage: instead of the moduli space $\M(\L)$ we consider the space of {\bf relative} realizations $C(\L,Z)$, see Definition \ref{2.7} for details.} 

Let $S\subset \R^m$ be a compact real-algebraic set, i.e. it is  
the zero set of a polynomial function $f: \R^m \to \R^n$. 
The set $S$ is contained in an $r$-ball $B_r(\O)$ centered at $\O$. We then 
apply Theorem A (the real case) and construct a functional linkage $\L$ for the pair 
$(f, B_r(\O))$. 
We let $\L_0$ be the abstract linkage obtained from $\L$ by gluing the output 
vertices to the base-vertex $v_1$. Let $p_0$ denote the restriction of  
the input mapping $p$ of $\L$ to $\M(\L_0)$, this map is the ``input map'' of 
$\L_0$. We show in \S \ref{Kemp} that $p_0$ is  
an analytically trivial polynomial covering over $S$. We obtain 

\bigskip
\noindent {\bf Theorem B.} {\em Let $S$ be any compact real-algebraic 
subset of $\R^m$. Then there is a linkage $\L_0$ so that $\M(\L_0)$ is 
Nash isomorphic to a disjoint union of a finite number of copies of $S$.} 

\medskip
\begin{rem}
Nash isomorphism is defined in \S 2. Nash isomorphism implies real analytic 
isomorphism. 
%In fact we construct $\L_0$ such that the input map 
%$p_0$ is an analytically trivial polynomial covering 
%$p_0: \M(\L_0)\to S$.
\end{rem}

Similarly we have 

\bigskip
\noindent {\bf Theorem B'.} {\em Let $S$ be any complex-algebraic 
subset of $\C^m$ and $U$ be an open (in the classical topology) bounded   
subset of $S$. Then there is a linkage $\L_0$ so that 
the input map $p_0$ is an analytically trivial polynomial covering over $U$.}

\begin{rem}
Here we treat $S$ as a {\bf real} algebraic set. 
\end{rem}

Now let $M$ be a compact smooth manifold. By work of Seifert, Nash, 
Palais and Tognoli (see \cite{AK} and \cite{T}), $M$ is diffeomorphic 
to a real algebraic set $S$, hence as a corollary of Theorem B we get 

\medskip
\noindent {\bf Corollary C.} {\em Let $M$ be a smooth compact manifold. Then 
there is a linkage $\L_0$ whose moduli space is diffeomorphic to a disjoint 
union of a number of copies of $M$.}

\begin{rem}
It is not true that for any compact smooth manifold $M$ there exists $\L_0$ 
such that $\M(\L_0)$ is diffeomorphic to $M$. This is because $\M(\L_0)$ 
admits a $\Z/2$-action (coming from $O(2)/SO(2)$) which is nontrivial 
provided that $\L_0$ is connected and $\M(\L_0)$ is not a point. Hence, if 
$M$ is compact, distinct from a point and does not admit a nontrivial 
$\Z/2$-action then $M$ cannot be diffeomorphic to the moduli space of a 
planar linkage. 
\end{rem}

Universality theorems similar  to Theorems A, B and Corollary C hold for 
moduli spaces of realizations of abstract arrangements in $\P^2$, see 
\cite{Mnev}, \cite{KM6} and \S \ref{arrangements} of this paper for 
more details. However we proved stronger realizability theorems in this case. In particular we get functional arrangements such that the input map is 
{\em injective}. 
In \S \ref{88} we establish a relation between the two kinds 
of universality theorems using moduli spaces of spatial Euclidean 
linkages, namely we show that the universality theorem for arrangements 
(Theorem \ref{F}) implies the following 

\medskip
{\bf Theorem D.} 
{\em Let $S$ be a compact real algebraic set defined over $\Z$. Then there are 
abstract linkages $\L, \L'$ so that: 

\smallskip
(1) $\M_0(\L, \R\P^2)$ is entire birationally isomorphic to $S$. 

(2) $\M_0(\L', \R^3)$ is an analytically trivial entire rational 
covering of $S$. 

\smallskip
\noindent Both $\M_0(\L, \R\P^2)$, $\M_0(\L', \R^3)$ are Zariski open and 
closed subsets in the moduli spaces $\M(\L, \R\P^2)$, $\M(\L', \R^3)$ of 
realizations of $\L, \L'$ in $\R\P^2$ and $\R^3$.}

\begin{rem}
It seems surprising that one can prove a somewhat stronger realization 
theorem for arrangements than for linkages. One explanation for this is 
that the image of the input map of any connected functional linkage is bounded. 
By a theorem of Sullivan \cite{Sul} a manifold with nonempty boundary 
can not be an algebraic set\, \footnote{The Euler characteristic of the 
link of a boundary point is $1$.}. Thus (unlike the case of functional 
arrangements) {\bf there are no functional 
linkages if we require the input map to be injective}.    
\end{rem}

\medskip
There is a long history of previous work on mechanical linkages 
(see \S \ref{hist} for a more detailed discussion). 
In particular, versions of Theorem A (for polynomial 
functions $\R \to \R^2$) and of Theorem E below   
were first formulated by A.~B.~Kempe in 1875 \cite{Kempe}, 
however, his statement of the theorem was vague and 
as far as we can tell, his proof requires 
corrections (due to possible degenerate configurations). 

\medskip
\noindent {\bf Theorem E.} (See \S \ref{draw}.) 
{\em Let $f= f(z, \bar{z}), f:\C \to \R$ be a 
polynomial function of the variables $z, \bar{z}$ and 
$\Ga:= f^{-1}(0)\subset \C$ be a real-algebraic curve. Pick an open 
(in the classical topology) bounded subset $U\subset \Ga$. Then there 
is a closed $\C$-functional linkage $\L_0$ so that the input map  
$p_0: C(\L_0, Z)\to \C$ is an analytically trivial polynomial covering 
over $U$. Thus we can ``draw'' arbitrary algebraic curves in $\R^2$ using 
planar linkages.} 

\medskip
The main problem with Kempe's proof is that it works well only for a 
certain subset of the moduli space, however near certain ``degenerate'' 
configurations the moduli space splits into several components and 
the linkage fails to describe the desired 
polynomial function (see for instance \S \ref{notion}). We use the 
``rigidified parallelograms'' to resolve this problem and get rid 
of the undesirable components.  

Kempe's methods were also insufficient to prove Theorem B and Corollary C 
even if the problem of ``degenerate configurations'' 
is somehow resolved. The second obstacle in  proving Theorem B is that the restriction $p_0$ of the regular ramified covering 
$p: \M(\L)\to Dom(\L)$ to  $\M(\L_0)$ apriori 
does not have to be an analytically trivial covering:

(a) It is possible that $\M(\L_0)$ intersects the ramification locus of $p$;

(b) Even if $p_0$ is a {\em topologically} trivial covering it 
can fail to be {\em analytically trivial} covering because of 
``quasiwalls'' (see \S \ref{elementary} for various examples): 
think of the function  $x^3:\R\to \R$. 

Both problems of degenerate configurations and reflection symmetries 
of linkages were neglected (or incorrectly resolved) in the previous 
work we have seen (e.g. \cite{HJW} and \cite{JS}). 
The first precise formulation of a theorem of 
the above type was given by W.~Thurston--  who stated a version of 
Corollary C about 20 years ago and has given lectures on it since. 
He realized that such a theorem would follow by combining the 19-th 
century work on linkages (i.e. Kempe's theorem) with the work of Seifert, 
Nash, Palais and Tognoli. However, Thurston did not write up a proof so 
we have no way of knowing whether he overcame the problems discussed 
above in the 19-th century work on linkages.  
There is also ambiguity concerning which theorem Thurston formulated in 
his lectures, we heard three different versions from three sources. 
According to the most recent (April 1997) oral com\-mu\-ni\-ca\-ti\-on 
from Thurston, he can also prove Corollary C.

We would like to thank a number of people who helped us with this work. 
The authors thank H.~King, S.~Lillywhite and R.~Schwartz  
for helpful conversations about real algebraic geometry and linkages. 
We are also grateful to M.~Karel, R.~Connelly and W.~Whiteley for 
supplying us with some of the references. The first author was supported by 
NSF grant DMS-96-26633 at University of Utah, the second author by NSF 
grant DMS-95-04193 at University of Maryland.  

\tableofcontents

\def\mf{\mathfrak }

\section{Some real algebraic geometry}

We begin by defining the categories of affine schemes and algebraic sets 
and the forgetful functor $\Phi$ from affine schemes to affine algebraic sets.

An {\em affine subscheme} ${\mathfrak X}$ of $\R^n$ is a locally ringed 
space (i.e. a topological space equipped with a sheaf of local rings) of 
the form $Spec\R[X_1,..., X_n]/I$ where $I$ is an ideal (see 
\cite[page 72]{Ha} for the definition of $Spec$). An {\em affine scheme} 
defined over $\R$ is a locally ringed space of the above form for some $n$. 

Let ${\mf X}=Spec\R[X_1,...,X_n]/I$ and ${\mf Y}=Spec\R[Y_1,...,Y_m]/J$. Then 
according to \cite{Ha} a morphism ${\mf f}: {\mf X}\to {\mf Y}$ 
consists of a pair $(f, \t f)$ where $f: {\mf X}\to {\mf Y}$ is a map of sets and 
$\t f: \O_{\mf Y}\to f_*\O_{\mf X}$ is a map of sheaves. However for the 
case in hand, $\t f$ is determined by the associated map of global sections
$$
\t f: \R[Y_1,...,Y_m]/J \to \R[X_1,...,X_n]/I \quad .
$$ 
Thus we may identify $\t f$ with a morphism $\t f: \R[Y_1,...,Y_m]\to 
\R[X_1,...,X_n]$ with $\t f(J)\subset I$. 

An affine scheme as above is said to be {\em reduced} if $I$ is equal 
to its radical $\sqrt{I}$. Recall that 
$$
\sqrt{I}= \{ g\in \R[X_1,..., X_n]: g^k \in I \hbox{~~for some~~} k\}
$$
An affine scheme ${\mf X}$ over $\R$ is said to be {\em  real reduced} 
if it is reduced and moreover $I$ as above has the property:

Suppose $g_1,..., g_{\ell}\in \R[X_1,...,X_n]$ satisfy 
$g_1^2 +...+g_{\ell}^2\in I$. Then $g_1,..., g_{\ell}\in I$. 

\medskip  
\noindent Let $I\subset \R[X_1,...,X_n]$ be an ideal. We define the 
{\em real radical} 
$\sqrt[r]{I}$ by
$$
g, g_1,..., g_{\ell} \in \sqrt[r]{I} \iff 
g^{2k}+ g_1^2 +...+g_{\ell}^2\in I \hbox{~~for some~~} k
$$
Note that $I\subset \sqrt{I}\subset \sqrt[r]{I}$. Thus 
${\mf X}= Spec\R[X_1,...,X_n]/I$ is real reduced if and only if $I= \i$. 
An ideal $I$ is {\em real radically closed} if $I=\i$. 

We now define the category of real algebraic sets. Let $I$ be an ideal in 
$\R[X_1,...,X_n]$. We define a subset $Z(I)\subset \R^n$ by
$$
Z(I)= \{x\in\R^n : g(x)=0, g\in I\}
$$
We note that $Z(I)= Z(\i)$. A subset $X\subset \R^n$ is said to be an 
{\em algebraic subset} if there exists $I\subset \R[X_1,...,X_n]$ with 
$X= Z(I)$. A set $X$ is said to be a real algebraic set if it is an 
algebraic subset of $\R^n$ for some $n$. Let $X$ and $Y$ be algebraic sets. 
Then a map $f:X\to Y$ is a {\em morphism} if there exist embeddings as 
above $X\subset \R^n, Y\subset \R^m$ such that $f$ is the restriction 
of a polynomial mapping $\t f: \R^n \to \R^m$. 

Similarly we define {\em semi-algebraic subsets} of $\R^n$. The collection 
of semi-algebraic subsets in $\R^n$ is the boolean algebra containing 
all sets of the form $\{x\in \R^n| f(x) >0\}$ for arbitrary polynomial 
functions $f$. For instance, any algebraic set $\{x| f(x)=0\}$ 
(where $f$ is a polynomial) is semi-algebraic since it is the complement 
of $\{x| f(x)>0\} \cup \{x| f(x)< 0\}$. A morphism between 
the semi-algebraic subsets $X\subset \R^n, Y\subset \R^m$ is a map 
$f: X\to Y$ which is the restriction of a polynomial mapping 
$\t{f}: \R^n \to \R^m$.  

\medskip
Let $X\subset \R^n$ be any subset. We define an ideal $\I= \I(X) 
\subset \R[X_1,...,X_n]$ by
$$
\I(X):= \{g\in \R[X_1,...,X_n]: g(x)=0 \hbox{~~for all~~} x\in X\}
$$
If $X= Z(I)$ then $\I= \i$. This is the real Nullstellensatz for 
polynomials \cite{BE}.  

\begin{thm}
$\I (Z(I))= \i$. 
\end{thm}

\begin{cor}
$Z$ and $\I$ give an order reversing bijection between real 
radically closed ideals in $\R[X_1,...,X_n]$ and algebraic subsets in $\R^n$.  
\end{cor} 

We define a functor $\Phi$ from the category of real affine subschemes 
of $\R^n$ to algebraic subsets of $\R^n$ by $\Phi({\mf X})= X$ 
with $X= Z(I)$, where ${\mf X}= Spec\R[X_1,...,X_n]/I$. If 
${\mf f}: {\mf X}\to {\mf Y}$ is a morphism then 
$f=\Phi({\mf f}):X\to Y$ is the associated map of sets. 
We have a right inverse $\Psi$ to $\Phi$. Let $X$ be an algebraic 
subset of $\R^n$. Define a real reduced scheme ${\mf X}_{can}= \Psi(X)$ by
$$
{\mf X}_{can}= Spec\R[X_1,...,X_n]/\I(X)
$$
If $f:X\to Y$ is a morphism we define 
${\mf f}: {\mf X}_{can}\to {\mf Y}_{can}$ by ${\mf f}= (f, f^*)$ 
where $f^*$ is the pull-back map on functions. 

\begin{rem}
As noted in the introduction if $\L$ is an abstract linkage then 
there is a canonical affine subscheme ${\mf M}(\L)$ of $\R^{2r-4}$ such that
$$
\Phi({\mf M}(\L))= \M(\L)
$$
However ${\mf M}(\L)$ is not necessarily reduced or real reduced. 
So general polynomial map $f: \M(\L_1)\to \M(\L_2)$ will not 
necessarily be of the form $\Phi({\mf f})$ for a morphism 
${\mf f}: {\mf M}(\L_1)\to {\mf M}(\L_2)$. 
\end{rem}

Suppose that $X$ and $Y$ are algebraic sets and that $f: X\to Y$ is a 
morphism. Suppose further that we have chosen affine schemes 
${\mathfrak X}, {\mathfrak Y}$ with $\Phi({\mathfrak X})=X$ and 
$\Phi({\mathfrak Y})=Y$. We will say $f$ is a {\em scheme-theoretic} 
morphism if there is a morphism ${\mathfrak f}: {\mathfrak X}\to {\mathfrak Y}$ 
such that $\Phi({\mathfrak f})= f$.  

The next definitions follow \cite{AK}. Let $X\subset \R^n$ be an algebraic set. 
We define an {\em entire rational} function $f: X\to \R$ to be a function 
which is locally (in a Zariski open \nbd of each point of $X$) the quotient 
of polynomials. Now let $X\subset \R^n, Y\subset \R^m$ be algebraic subsets. 
An {\em entire rational map} $X\to Y$ is a mapping of sets where components are 
entire rational functions. 
A {\em entire birational isomorphism} $f: X\to Y$ is an entire 
rational map which has entire rational inverse (in particular $f$ is a 
homeomorphism). Note that the notion of {\em entire birational isomorphism} 
is more restrictive than birational isomorphism 
(a birational map does not have to be defined everywhere). 

In what follows we will need the notion of the fiber product of 
affine schemes defined over $\R$. Suppose that we have a diagram of real 
affine schemes 
\[
\begin{array}{ccc}
~ & ~ & {\mf Y}\\
~ & ~ & {\mf g}\downarrow\\
{\mf X} & \stackrel{{\mf f}}{\longrightarrow} & {\mf Z}\\
\end{array}
\]
The categorical fiber product ${\mf X}\times_{\mf Z}{\mf Y}$ is an 
affine scheme ${\mf W}$ defined over $\R$ with morphisms 
$\pi_1:{\mf W}\to {\mf X}$ and $\pi_2: {\mf W}\to {\mf Y}$ such 
that we have a diagram
\[
\begin{array}{ccc}
{\mf W} & \stackrel{\pi_2}{\lra} & {\mf Y}\\
\pi_1\downarrow & ~ & {\mf g}\downarrow\\
{\mf X} & \stackrel{{\mf f}}{\longrightarrow} & {\mf Z}\\
\end{array}
\]
satisfying the universal property that:

For any affine scheme ${\mf V}$ defined over $\R$ the natural map of sets
$$
Mor({\mf V}, {\mf W})\to Mor({\mf V}, {\mf X}) \times Mor({\mf V}, {\mf Y})
$$
is an injection with image the subset of pairs $(\al,\be)$ 
satisfying ${\mf f}\circ \al= {\mf g}\circ \be$. 

We now have 

\begin{lem}
The categorical fiber product ${\mf X}\times_{\mf Z} {\mf Y}$ 
exists (as an affine scheme defined over $\R$) and is unique up to 
canonical isomorphism. 
\end{lem}
\proof Choose representations 
$$
{\mf X}= Spec\R[X_1,...,X_n]/I, {\mf Y}= Spec\R[Y_1,...,Y_m]/J, 
{\mf Z}= Spec\R[Z_1,...,Z_{\ell}]/K
$$
Then define 
$$
{\mf X}\times_{\mf Z} {\mf Y}:= Spec(\R[X_1,...,X_n]/I 
\otimes_{\R[Z_1,...,Z_{\ell}]/K} \R[Y_1,...,Y_m]/J)
$$
This proves existence. Uniqueness is obvious.  \qed  

\begin{rem}
In down-to-earth terms ${\mf X}\times_{\mf Z} {\mf Y}$ 
is represented by a subscheme of $\R^{n+m}$ 
(with coordinates $(X_1,...,X_n, Y_1,...,Y_m)$) 
defined by the union of equations defining ${\mf X}$ and ${\mf Y}$ 
together with the equations
$$
f^*Z_i = g^* Z_i, i=1,...,\ell
$$
\end{rem}

It is remarkable that the abstract notion of fiber product of schemes 
corresponds to the operation of gluing of linkages, see Section \ref{fun}. 

We next define a functor $\Theta$ from the category of affine schemes 
over $\R$ to real analytic spaces. The definition is complicated because we 
do not want to assume that our real analytic spaces are reduced- we do not want 
to loose track of nilpotents. 

Let ${\mathfrak X}=Spec(\R[X_1,..., X_n]/I)$ be a subscheme of 
$\R^n$. Our goal is to define the sheaf $\O^{an}$ (in the classical 
topology) of real analytic ``functions'' (they may be nilpotent) over 
the space of real maximal ideals ${\mathfrak X}_{max}$ in the ring 
$\R[X_1,..., X_n]/I$. We first define $\O^{can}$ for 
the affine space $\R^n$. Let $U\subset \R^n$ be an open subset. 
Then $\O^{an}(U)$ is defined to be the ring of all analytic functions 
from $U$ to $\R$. Now the subscheme ${\mf X}\subset \R^n$ determines 
an ideal sheaf $\J$ of $\O_{\R^n}|{\mf X}$   such that we have an 
exact sequence of sheaves over ${\mf X}$:
$$
0\to \J \to \O_{\R^n}|{\mf X} \to \O_{{\mf X}} \to 0
$$
(note that the restriction is the ``abstract'' restriction, 
\cite[page 65]{Ha}). The sheaf of ideals $\J \subset \O_{\R^n}|{\mf X}$ 
induces a sheaf of ideals $\O_{\R^n}|{\mf X}_{max}$ (see \cite[page 78]{Ha}) 
which extends to a sheaf of ideals 
$\J^{an}\subset \O_{\R^n}^{an}|{\mf X}_{max}$. We define 
$$
\O^{an}_{{\mf X}}:= (\O_{\R^n}^{an}|{\mf X}_{max})/\J^{an}
$$

\begin{rem}
This definition states that an analytic function $f\in \O_{\R^n}^{an}$ 
``vanishes''  on $X\cap U$ if it may be written as $f= \sum_{i=1}^N g_i f_i$ 
where $I=(f_1,...,f_N)$ and $g_1,...,g_N$ are analytic functions on $U$. This 
is a stronger requirement than requiring that the induced {\bf function} $f$ 
on $X\cap U$ is identically zero.   
\end{rem}

In the case that ${\mf X}={\mf X}_{can}$ is real reduced the definition of 
$\O^{an}$ is 
much simpler. Let $x\in {\mf X}$. Then  $f\in \O^{an}_{{\mf X},x}$ if 
there is a \nbd $U$ of $x$ in $\R^n$ (in the classical topology) and an 
analytic function $\t{f}$ on $U$ such that $\t{f}|{\mf X}$ agrees with $f$.  
If $X$ is an algebraic subset 
of $\R^n$ we obtain the sheaf of analytic functions $\O^{an}$ on $X$ by 
considering the real  reduced scheme ${\mf X}_{can}$ and using the 
definition at the beginning of this paragraph. 

We obtain the corresponding notions of analytic morphism. Let $X\subset \R^n$ and 
$Y\subset \R^m$  be algebraic subsets. Then an {\em analytic morphism}
$f: X\to Y$ is a map such that for each $x\in X$ there is an open \nbd (in the 
classical topology) $U$ of $x$ in $\R^n$ such that $f|U\cap X$ is the restriction 
of a real analytic map $\t{f}: U\to \R^m$. An {\em analytic isomorphism} is a morphism 
which has an analytic inverse. 

More generally, assume we are given affine subschemes ${\mf X}\subset \R^n$, 
${\mf Y}\subset \R^m$ with $\Phi({\mf X})= X, \Phi({\mf Y})=Y$. We define a morphism
$$
{\mf f}: {\mf X}^{an}\to {\mf Y}^{an}
$$
to be a morphism of locally-ringed spaces \cite[page 72]{Ha}. Such a morphism 
${\mf f}$ induces a morphism ${\mf f}_{can}$ of real reduced spaces and 
consequently determines an analytic morphism $f:X\to Y$.   We say that an 
analytic morphism $f: X\to Y$ is a {\em scheme-theoretic} analytic morphism 
if it is induced by ${\mf f}: {\mf X}^{an}\to {\mf Y}^{an}$ as above. 
We define a {\em scheme-theoretic} analytic isomorphism to be a 
 {\em scheme-theoretic} analytic morphism which has a 
{\em scheme-theoretic} analytic inverse. 

In what follows it is important to note 

\begin{lem}
A morphism 
$$
{\mf f}: {\mf X}^{an}\to {\mf Y}^{an}
$$
induces a morphism 
$$
{\mf f}_{can}: {\mf X}^{an}_{can}\to {\mf Y}^{an}_{can}
$$
\end{lem}
\proof Indeed, $f$ induces  a map of point sets $f: X\to Y$. Also, 
given an open subset $U\subset Y$, $f$ comes with a local homomorphism
$$
\t{f}_U: \Ga(U, \O^{an}_{\R^m}|U) \to \Ga(f^{-1}(U), \O^{an}_{\R^n}|f^{-1}(U)) 
$$  
which commutes with point evaluation since $\t{f}_U$ is local. 
Hence $\t{f}_U$ carries $\Ga(U, \J_{can,Y})$ into $\Ga(f^{-1}(U), \J_{can,X})$ 
and we obtain an induced map of quotients.  \qed  

\begin{rem}
The previous  lemma also follows from the local Nullstellensatz 
\cite[Proposition 2.8 (h), page 216]{ABR} which implies 
$$
\J_{can}= \sqrt[r]{\J}
$$
where $\J_{can}\subset \O^{an}_{\R^n}|X$ is the ideal sheaf of 
${\mf X}^{an}_{can}$ (i.e. the sheaf of analytic functions vanishing on 
$X$ in the usual sense). The above equation means that the maximal 
ideals $\J_{can,x}$ and $\J_x$ in $\O^{an}_{\R^n,x}$ are related by 
$$
\J_{can,x}= \sqrt[r]{\J_x}
$$
\end{rem}

\begin{lem}
\label{locis}
Assume we are given affine schemes ${\mf X}$ and ${\mf Y}$ with $\Phi({\mf X})=X$ 
and $\Phi({\mf Y})=Y$. Suppose that $f: X\to Y$ is a bijective 
scheme-theoretic analytic morphism which induces isomorphisms of the 
completed local $\R$-algebras:
$$\widehat{O_{{\mathfrak Y},y}}\to \widehat{O_{{\mathfrak X},x}} $$ 
for each $x\in X$ with $y=f(x)\in Y$. Then $f$ is a scheme-theoretic 
analytic isomorphism.   
\end{lem}
\proof The lemma follows from Artin's theorem, see \cite[Theorem 3.1]{GM}.  
\qed

Let ${\mf X}$ be an affine scheme over $\R$ and $x\in X$. Let ${\mf m}\in 
\widehat{\O_{{\mf X},x}}$ be a maximal ideal. We define the $N$-th order 
tangent space $T^{(N)}_x({\mf X})$ to be the set
$$
Hom_{\R-alg}(\widehat{\O_{{\mf X},x}}, \R[t]/(t^{N+1}))
$$
A scheme-theoretic analytic morphism $f: X\to Y$ with $f(x)=y$ induces maps 
$$
(D^N f)_x: T^{(N)}_x{\mf X} \to T^{(N)}_y{\mf Y}
$$
We will use the following lemma to verify that a scheme-theoretic 
analytic morphism is a scheme-theoretic analytic isomorphism and hence 
induces an analytic isomorphism of germs
$$
f: (X,x)\to (Y,y)
$$ 

\begin{lem}
\label{algeo}
(See \cite[Lemma 2.1]{KM6}.) Suppose that 
${\mf f}: ({\mf X}, x)\to ({\mf Y},y)$ is a 
morphism of germs such that the maps  
$(D ^{(N)}{\mf f})_x$ are bijective for all $N$. Then 
${\mf f}$ (and consequently ${\mf f}_{can}$) is an isomorphism of germs. 
\end{lem}

The following lemma is useful in verifying the hypothesis of Lemma \ref{algeo}

\begin{lem} 
(See \cite[Remark 4.3]{GM2}.) Suppose ${\mf f}: ({\mf X},x)\to ({\mf Y},y)$ 
is a morphism of germs such that the induced map 
$D{\mf f}_x: T_x({\mf X})\to T_y({\mf Y})$ is injective. Then the 
maps $(D^N{\mf f})_x$ are injective for all $N$. 
\end{lem}

Let $f: X\to Y$ be a scheme-theoretic analytic morphism where 
$X= \Phi({\mf X})$, $Y= \Phi({\mf Y})$. An {\em irregular  
point} of $f$ is a point where $f$ is not a local analytic 
isomorphism (in the scheme-theoretic sense). 
%A point $y\in Y$ will be called {\em a regular value of} $f$ 
%if each $x\in f^{-1}(y)$ is not critical.  
To check that $x\in X$ is not irregular   
it is enough to verify that $f$ induces bijections of Zariski 
tangent spaces of all orders: $T^{(N)}_{x}{\mf X}\to T^{(N)}_{f(x)}{\mf Y}$, 
see Lemma \ref{algeo}.

\begin{defn}
\label{cov}
Suppose that $X, Y$ are real algebraic sets. Then a finite 
{\bf analytically trivial} covering $f: X\to Y$ is an analytic  
map such that the restriction of $f$ to each connected component of $X$ 
is an analytic isomorphism. 
  
We say that $f: X\to Y$ is an analytically trivial {\bf polynomial} 
covering if it is an polynomial morphism which is an analytically trivial 
regular covering whose group $G$ of deck transformations 
consists of polynomial automorphisms. We retain 
the name {\bf analytically trivial polynomial covering} for restriction of 
such $f$ to a $G$-invariant open subset\footnote{With respect to 
the classical topology.} of $X$.  
\end{defn}

Note, that we do not claim here that $X$ splits into disjoint union of 
Zariski components each of which is polynomially isomorphic to $Y$. 
It might  happen that the real-algebraic set $X$ is irreducible, but 
$f$ is not 1-1. 

We will now give a version of Definition \ref{cov} in terms of 
{\em Nash functions} (see Lemma \ref{nas}). Let $X$ be a real semi-algebraic 
set and $U\subset X$ be an open subset (in the classical topology). 

\begin{defn}
A function $f:U \to \R$ is {\bf Nash} if it is real-analytic and 
there exist polynomial functions $p_0, p_1,..., p_d$ 
not all equal to zero such that the equation
$$
p_0 + p_1 f + ... + p_d f^d=0
$$
holds identically on $U$. 
\end{defn}

We define the sheaf ${\cal N}_X$ on $X$ by defining ${\cal N}_X(U)$ 
to be the $\R$-algebra of Nash functions on $U$. Now let $X$ and $Y$ 
be real algebraic sets and $f: X\to Y$ be a continuous map. Define $f$ 
to be a {\em Nash morphism} if $f^* {\cal N}_Y \subset f_*{\cal N}_X$. 
An equivalent definition is the following. Choose an embedding 
$Y\subset \R^n$. Then $f$ is Nash if and only if the components of $f$ are 
Nash functions on $X$. We have the following useful criterion 
for an analytic function to be Nash:

\begin{lem}
\label{N1}
(See \cite[Proposition 8.1.7]{BCR}.) 
Suppose $X$ is a real semi-algebraic set and $f: X\to \R$ is an 
analytic function on $X$. Then $f$ is Nash if and only if the 
graph $Gr(f)$ of $f$ 
is a semi-algebraic subset of $X\times \R$. 
\end{lem}

We will also need

\begin{lem}
(See \cite[Theorem 2.4.5]{BCR}.) 
Suppose $X$ is an real semi-algebraic set. Then the topological 
components of $X$ are semi-algebraic sets. 
\end{lem}

Now we can prove the result we need:

\begin{lem}
\label{nas}
Suppose $X$ and $Y$ are real algebraic sets and $X= \cup_{i=1}^k X_i$ is 
the decomposition of $X$ into connected components. Suppose 
$f:X\to Y$ is a polynomial map such that each $f_i:= f|X_i$ 
is an analytic isomorphism. Then 
each $f_i$ is a Nash isomorphism.  
\end{lem}
\proof Choose embeddings $X\subset \R^m, Y\subset \R^n$. Then we have:
$$
Gr(f)= \cup_{i=1}^k Gr(f_i)
$$ 
where $Gr(f_i)$ is a semi-algebraic subset of $X\times Y$. Let 
$g_i:= f_i^{-1}$. Then $Gr(g_i)\subset Y\times X$ is the image of 
$Gr(f_i)$  under the map which exchanges $X$ and $Y$. Hence $Gr(g_i)$ 
is semi-algebraic. Let $\pi_j: \R^m \to \R$ be the $j$-th coordinate 
projection and $\Pi_j: \R^n\times \R^m \to \R^n\times \R$ be the 
projection defined by $\Pi_j(y,x)=(y, \pi_j(x))$. Then 
$Gr(\pi_j\circ g_i)= \Pi_j(Gr(g_i))$ is the image of a 
semi-algebraic set under a projection. Hence  $Gr(\pi_j\circ g_i)$ 
is semi-algebraic by \cite[Theorem 2.2.1]{BCR}. Therefore $\pi_j\circ g_i$ is 
Nash by Lemma \ref{N1} which implies that $g_i$ is a Nash morphism and 
hence isomorphism. \qed
  
\medskip
Thus, if $f: X\to Y$ is an analytically trivial polynomial covering then 
$X$ is Nash isomorphic to a disjoint union of a finite number of copies of $Y$. 

\medskip
We will identify $\R^n$ with the affine part of $\R\P^n$. Suppose that 
$X\subset \R^n$ is an affine real algebraic set. Then $X$ is said to be 
{\em projectively closed} it its Zariski closure in $\R\P^n$ equals $X$. 
Clearly each projectively closed subset must be compact (in the classical 
topology). It turns out that the converse is ``almost true'' as well:

\begin{thm}
\label{closed}
(Corollary 2.5.14 of \cite{AK}) 
Suppose that $X\subset \R^n$ is a compact affine algebraic set. 
Then $X$ admits an entire rational isomorphism 
to a projectively closed affine 
algebraic subset $X'$ of $\R^n$. 
Moreover, if $X$ is defined over $\Z$ then $X'$ is defined over  
$\Z$ as well. 
\end{thm}

\noindent We will need the following theorem  which is a modification of 
 \cite[Corollary 2.8.6]{AK} or \cite{T}:

\begin{thm}
\label{nash}
(Seifert-Nash-Palais-Tognoli) Suppose that $M$ is a smooth compact 
manifold~\footnote{Not necessarily connected.}. Then $M$ is diffeomorphic 
to a projectively closed real affine algebraic set $S$.% defined over $\Z$. 
\end{thm}

\begin{figure}[tbh]
\leavevmode
\centerline{\epsfxsize=3in\epsfbox{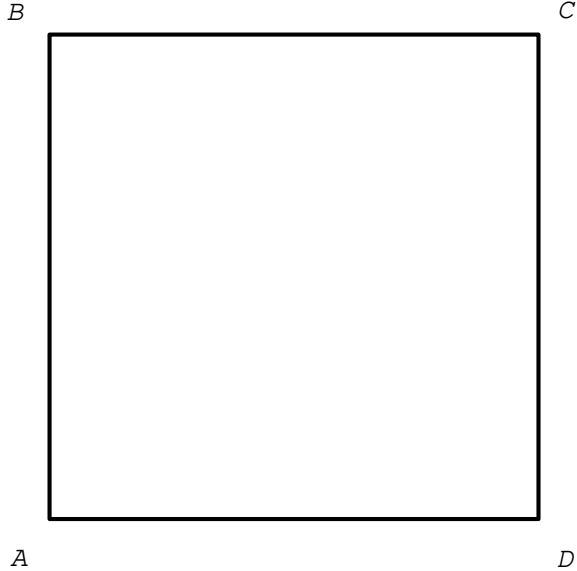}}
\caption{\sl The abstract square (\`a la Malevich).}
\label{Fig1}
\end{figure}

\begin{rem}
This theorem is a combination of \cite[Corollary 2.8.6]{AK} and 
Theorem \ref{closed} (for the assertion that $S$ is projectively closed). 
%We are grateful to H.~King 
%and defined over $\Z$. 
\end{rem}

\begin{Not}
If $f: X\to X$ then $Fix(f)$ will denote the fixed point set of $f$. 
\end{Not}

\begin{Not}
Let $z\in \k^n$, $r>0$. Then $B_r(z)$ will denote the disk in $\k^n$ 
of the radius $r$ and center at $z$. 
\end{Not}

A {\em domain} in $\k^N$ ($k=\R$ or $\k= \C$) is a subset $V$ with nonempty interior, 
where we use the classical topology.

\section{Linkages}

\subsection{Abstract linkages and their realizations}
\label{notion}

Throughout this paper we shall use the notation $\V(L)$ for the set 
of vertices of a graph $L$ and $\E(L)$ for the set of edges of $L$.

\begin{defn}
\label{2.1}
An {\bf abstract marked linkage} $\L$ is a triple $(L, \ell, W)$ 
consisting of a %connected 
graph $L$, an ordered subset $W\subset \V(L)$ 
and a positive function $\ell: \E(L)\to \R_+$ (a metric on $L$). 
We  shall assume that no vertex of $L$ is connected by an edge to itself. 
The elements of $W$ are called the {\bf fixed vertices} of $\L$ and 
the choice of $W$ is called {\bf marking}. If $W$ is empty then we 
call $\L$ an {\bf abstract linkage}. 
A special case of abstract marked linkage is an abstract {\bf based} linkage 
where $W$ consists of two vertices $v_1, v_2$ connected by an edge $e^*$. 

\end{defn}

We do not require $\ell$ to define a {\em metric} on $L$: for instance 
the triangle with the edge-lengths $1, 1, 3$ satisfies our axioms.  
However in what follows we will refer to the pair $(L, \ell)$ as a 
{\em metric graph}. 

\begin{defn}
\label{2.2}
Let $\L= (L, \ell, W )$ be an abstract linkage. 
A {\bf planar realization}  of $\L$ is a map $\phi: \V(L)\to \R^2$ 
such that: if $v$ and $w$ are joined by an edge $[vw]$ then 
$$
|\phi(v) - \phi(w)|^2= (\ell[v w])^2
$$
We let $C(\L)= C(\L, \0)$ be the set of all planar realizations of $\L$, 
it is called the {\bf configuration space} of $\L$. 
%We shall use the notation 
%${\mf C}(\L)$ for the scheme corresponding to $C(\L)$. 
\end{defn}

\begin{defn}
\label{2.3}
Let $\L= (L, \ell, W )$ be an abstract linkage, $W= (v_1,..., v_n)$ be 
the marking. Let $Z= (z_1,..., z_n)\in \C^n$, called 
{\bf the image of marking}. A {\bf relative} {\bf planar realization} 
 of $\L$ is a realization  
$\phi\in C(\L)$ such that $\phi(v_j)= z_j$ for all $j$. 
We let $C(\L, Z)$ be the set of all relative planar realizations of $\L$, 
it is called the {\bf relative configuration space} of $\L$. 
\end{defn}

A special case of the relative configuration space is when $\L$ 
is a based linkage: 

\begin{defn}
\label{2.4}
Let $\L= (L, \ell, e )$ be a based linkage. We define the 
{\bf moduli space} $\M(\L)$ by 
$$
\M(\L)= \{\phi\in C(\L) : \phi(v_1)=(0,0),\ \phi(v_2)= (\ell(e), 0)\}
$$
\end{defn}

\begin{figure}[tbh]
\leavevmode
\centerline{\epsfxsize=3in\epsfbox{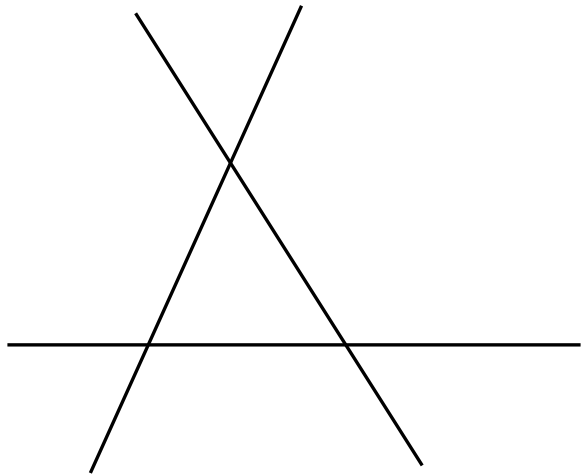}}
\caption{\sl The moduli space of the square.}
\label{Fig2a}
\end{figure}

If $L$ is connected then $\M(\L)$ is a compact real-algebraic subset 
of $\R^{2r-4}$, where $r$ is the number of vertices in $L$.  The 
algebraic set $C(\L)$  canonically splits as the product $\M(\L)\times 
E(2)$ (the group $E(2)$ of orientation-preserving isometries of $\R^2$ 
has obvious real-algebraic structure), thus we shall identify the quotient 
$C(\L)/E(2)$ and $\M(\L)$. Note that $\M(\L)$ admits an algebraic  
automorphism induced by the complex conjugation in $\C= \R^2$. 

As we shall see in Section \ref{fixing} for each $(\L, Z)$ 
there is a based abstract linkage $\a$ so that $\M(\a)$ is a 2-fold or 1-fold 
analytically trivial covering of $C(\L, Z)$.

Suppose that $\L'\subset \L$ is a sublinkage, i.e. 
$(L', \ell')\subset (L, \ell)$ is a subgraph such that $\ell'$ is the 
restriction of $\ell$ and the marking $W'$ of $\L'$ is the intersection 
$W\cap L'$. If $Z$ is the image of $W$, then we have naturally  
defined $Z'$, the image of $W'$. Thus we have natural 
scheme-theoretic {\em restriction} morphism 
$$
Res: C(\L, Z)\to C(\L', Z'), \quad \phi\mapsto \phi|L'
$$
In particular, if $\L'$ consists of a single vertex 
$v\in \L$ then we let $eval_v: C(\L, Z)\to \R^2$ be 
the evaluation map $\phi\mapsto \phi(v)$. We will identify $eval_v$ 
with the restriction map 
$$
Res: \phi \mapsto \phi|\{v\}\in C(\{v\})
$$

\begin{figure}[tbh]
\leavevmode
\centerline{\epsfxsize=4in\epsfbox{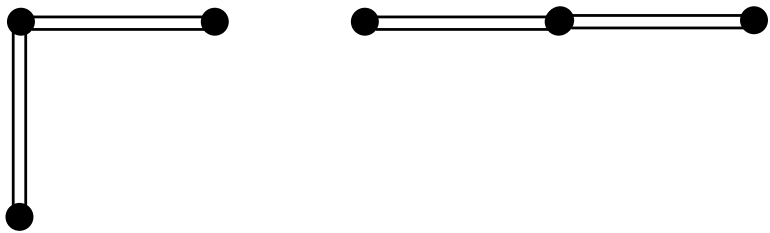}}
\caption{\sl Degenerate realizations of the abstract square. Black 
disks denote images of the vertices.}
\label{Fig3}
\end{figure}

Many of the problems with the 19-th century work on linkages can be 
traced to neglecting degenerate realizations of a square. 

\begin{defn}
A {\bf square} is the abstract polygonal linkage with four edges 
 all of which have equal length (see Figure \ref{Fig1}). We let $e^*:= [AB]$. 
More generally, an {\bf abstract parallelogram} $S$ is the polygonal 
linkage with four edges where the alternating edges have equal lengths. 
\end{defn}

We have

\begin{lem}
\label{2.5}
The moduli space of the square is isomorphic to a union 
of {\bf three} smooth curves of degree 2 (each one necessarily rational) in $\R\P^6$ such that each pair intersects in a point and at each point of intersection the tangent spaces have 2-dimensional span (see Figure \ref{Fig2a}). 
\end{lem}
\proof See \cite[Case III, page 120]{GN}. 

\medskip

\begin{figure}[tbh]
\leavevmode
\centerline{\epsfxsize=3in\epsfbox{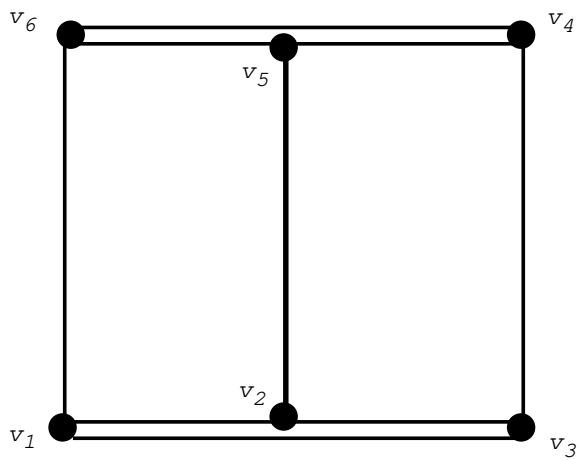}}
\caption{\sl The rigidified parallelogram $\Si$. 
We choose: $\ell[v_1 v_2]=\ell[v_2 v_3]= 
\ell[v_1v_3]/2= \ell[v_6 v_5]=\ell[v_5 v_4]=\ell[v_6v_4]/2$.}
\label{figs}
\end{figure}

Two of the components of the moduli space of the square consist of 
``degenerate'' realizations of the square (Figure \ref{Fig3}). 
We can eliminate the components consisting of degenerate squares 
by ``rigidifying'' the square as on Figure \ref{figs}. We 
rigidify parallelogram linkages in an analogous way. Henceforth 
all parallelogram sublinkages that appear in {\em elementary linkages} 
(Section \ref{elementary}) will be rigidified-- but we will not draw 
the extra edges. Notice that each rigidified parallelogram $\L$ contains 
two ``degenerate'' triangular linkages: one with 
the vertices $v_1, v_2, v_3$ and second with the vertices $v_4, v_5, v_6$. 
Thus the ring of the configuration scheme of the rigidified parallelogram 
has nilpotent elements. Since we are interested mostly in the reduced schemes 
we make the following  

\begin{figure}[tbh]
\leavevmode
\centerline{\epsfxsize=3in\epsfbox{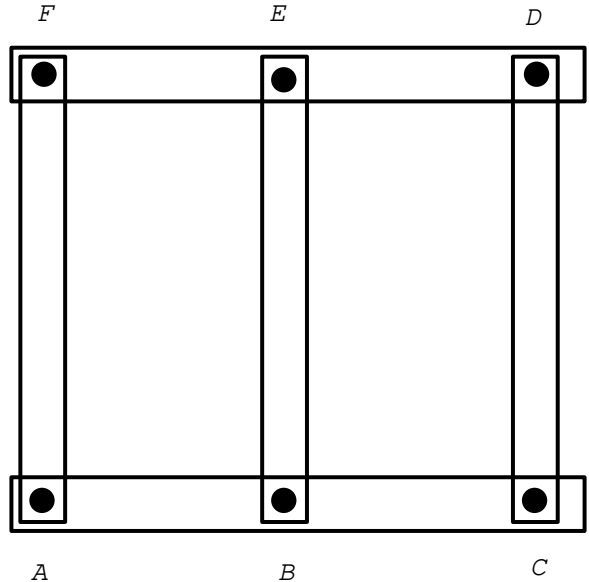}}
\caption{\sl A mechanical model for the rigidified parallelogram $\Si$.}
\label{mech}
\end{figure}

\begin{conv}
\label{contri}
Suppose that $\L$ is one of the {\bf elementary linkages} 
(translator, pantograph and inversor) in 
\S \ref{elementary} and $\Del\subset \L$ is a 
{\bf degenerate triangle}, i.e. a triangle with the vertices 
$A, B, C$ so that 
$$
\ell[AC]= \ell[AB]+ \ell[BC]
$$
Let $r:=\ell[AB]/\ell[AC], s:= \ell[BC]/\ell[AC]$. 
(It is clear that $r=s= 1/2$ for parallelograms 
but the construction is more general, the convention 
is used for the pantographs as well.) 
Then throughout the paper we will use the following scheme 
${\mathfrak C}(\L, Z)$ whose set of real points is the configuration 
space $C(\L, Z)$ of $\L$: we use in ${\mathfrak C}(\L, Z)$ the equations 
$$
\phi(B)= r\phi(A)+ s\phi(C)
$$
for each degenerate triangle instead of the equations
$$
|\phi(A)- \phi(B)|^2= \ell[AB]^2, \quad |\phi(C)- \phi(B)|^2= \ell[BC]^2
$$
This choice of the scheme is determined by the ``mechanical'' reasons: to 
make actual mechanical model of the abstract linkage $\L$ drill the 
hole $B$ in the bar $[AC]$ within the distances $\ell[AB]$ and $\ell[BC]$ 
from the holes $A$ and $C$ respectively. See Figure \ref{mech} for 
mechanical model of the rigidified parallelogram. We will use the notation 
${\mf M}'(\Si)$ for the moduli scheme of the rigidified parallelogram $\Si$ with 
the partially reduced structure as above. 
\end{conv}

Recall that $S$ is the (unrigidified) parallelogram linkage with 
$e^*=[AB]$. Then we have an embedding of affine {\em schemes} 
$i: \M'(\Si) \to \M(S)$. 

We let $\phi_1$ and $\phi_2$ be the degenerate realizations of the 
rigidified parallelogram $\Si$ (which could be the rigidified square). 
The following lemma will be very important in what follows.

\begin{lem}
\label{2.6}
The real reduced structure on ${\mf M}'(\Si)$ is a projective line where 
real points correspond to (convex) parallelograms. ${\mf M}'(\Si)$ has 
exactly two singular real points, the degenerate parallelograms 
$\phi_1, \phi_2$.  
\end{lem}
\proof Recall that the distance function between two oriented 
straight lines in the Euclidean plane is convex and it is strictly 
convex unless these lines are parallel. Thus, if $A\ne B$ and $C\ne D$ are 
points in $\C$ such that 
$$
\|A-D\|= \|B-C\|= \|(A+B)/2 - (D+C)/2\|
$$
then the lines through $A, B$ and $D,C$ are parallel and 
$$
A- B = D- C 
$$
Therefore real points $\phi$ of ${\mf M}'(\Si)$ correspond to parallelograms:
$$
\phi(v_1)- \phi(v_3)= \phi(v_6) - \phi(v_4)
$$
The assertion that the only singular real points of ${\mf M}'(\Si)$ are 
degenerate realizations can be proven analogously to \cite{KM2}. 
An alternative proof follows from the discussion below. 

\begin{figure}[tbh]
\leavevmode
\centerline{\epsfxsize=2.5in\epsfbox{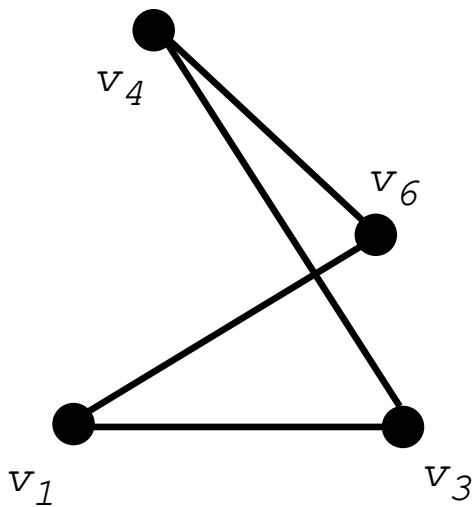}}
\caption{\sl An antiparallelogram.}
\label{fig:6}
\end{figure}

We first consider the case of a parallelogram $S$ which is not a square. 
The moduli space $\M(S)$ for such $S$ is described in 
\cite[Case II, page 120]{GN}. The authors of \cite{GN} 
describe the (projectivized) moduli space as a real projective 
subvariety of ${\mathbb P}^6$. They find that the moduli space 
is the union of a smooth curve of degree two (necessarily isomorphic 
to ${\mathbb P}^1$) and a smooth curve of degree four which is also 
isomorphic to  ${\mathbb P}^1$ (see \cite[page 119]{GN} where 
it is proved that if the degree is four then the genus is zero). 
Moreover they 
prove that the real points of the quartic correspond to the set of 
``antiparallelogram'' (see Figure \ref{fig:6}) 
realizations of the linkage. The authors also 
prove that the components of $\M(S)$ intersect in two points, the 
two degenerate parallelograms. 
It now follows from the paragraph above that $i(\M'(S))$ is the 
quadratic curve $C$. We leave it to the reader to verify that the points 
$\phi_1,\phi_2$ are singular 
on $\M'(\Si)$. The lemma follows in the parallelogram case.   

We leave the proof of the lemma in the case of a square to the reader, 
it is a consequence of \cite[case III, page 120]{GN}. In this case 
$\M(S)$ is the union of three quadratic curves, two of the components 
correspond to degenerate realizations, see Figures \ref{Fig2a}, \ref{Fig3}.  The image of $i$ is the component consisting of rhombi. \qed 

The (relative) configuration space space of a linkage $\L$ is 
naturally the set of real points of an affine subvariety in $\A^{2r}$ 
where $r$ is the number of vertices in $L$. Hence for each realization 
$\phi\in C(\L, Z)$ the Zariski tangent 
space\footnote{It and higher  order Zariski tangent spaces depend on the 
choice of the affine scheme, see 
the Convention \ref{contri}.} 
$T_{\phi}{\mf C}(\L, Z)$ embeds into $\R^{2r}$. We will refer to 
elements of  $T_{\phi}{\mf C}(\L, Z)$ as {\em infinitesimal deformations} of 
$\phi$. If $w$ is a vertex of $L$ and 
$\xi$ is an infinitesimal deformation of $\phi$ then we get a vector 
$\xi|w\in \R^2$ given by projection of $\xi$ to the appropriate 
$\R^2$-factor of $\R^{2r}$.   We use similar definitions for higher 
order infinitesimal deformations, i.e. elements of the $n$-th order Zariski 
tangent space $T^{(n)}_{\phi} {\mf C}(\L, Z)$. 
If $\xi\in T^{(n)}_{\phi}{\mf C}(\L, Z)$ 
then $\xi|w$ is the image of $\xi$  under the push-forward map 
$$
(eval_v)_* : T^{(n)}_{\phi} {\mf C}(\L, Z) \to T^{(n)}_{\phi} {\mf C}(\{w\})
$$

The proof of the following lemma is straightforward and is left to the reader:  

\begin{lem}
\label{tri}
Let $\L$ be a based triangular linkage with the vertices 
$v_1, v_2, v_3$, $e^*:= [v_1v_2]$ so that 
$\ell[v_1v_2]< \ell[v_2v_3]+ \ell[v_1v_3]$, 
$\ell[v_1v_3]< \ell[v_2v_3]+ \ell[v_1v_2]$
and $\ell[v_2v_3]< \ell[v_1v_2]+ \ell[v_1v_3]$ (i.e $\L$ is nondegenerate 
and the function $\ell$ defines a metric on $L$). 
Then $\M(\L)$ consists of two points; if 
$\xi\in T^{(n)}_{\phi}{\mf C}(\L, Z)$ is an infinitesimal deformation of 
order $n$ such that $\xi|v_1=\xi|v_2=0$, then $\xi=0$.  
\end{lem}

\begin{Not}
Throughout the paper we shall use the following notations:

\begin{itemize}
\item If $A, B$ are distinct points in the plane then $(AB)$ will denote the 
straight line through $A, B$.  
\item If $A, B, C$ are points on the plane then $\Del(A,B,C)$ will denote the 
triangle with the vertices $A,B,C$. 
\item If $(A, B, C, D)$ is a quadruple of points on the plane then $[ABCD]$ 
(sometimes we will also put commas between the letters) will  denote the 
quadrilateral with the given 
vertices (they are connected by edges according to the cyclic order). 
\end{itemize}
\end{Not}

\subsection{Fiber sums of linkages}

The operation of fiber sum of linkages is analogous to generalized 
free products of groups (i.e. amalgamated free product or HNN-extension). 
Let $\L'= (L', \ell', W')$, $\L''= (L'', \ell'' , W'')$ be abstract marked  
linkages. Suppose that we have a map $\beta$ between (nonempty) 
subsets of vertices 
$$
\beta: S'\subset \V(L') \to S''\subset \V(L'') 
$$
If the images $Z', Z''$ of $W', W''$ are given we require 
$$
\phi'(w_j)= \phi''(\be (w_j))
$$
for each $w_j\in W'$ and $\phi'\in C(\L', Z'), \phi''\in C(\L'', Z'')$. 

Then the {\em fiber sum} $\L$ of linkages $\L', \L''$ associated 
with $\beta$ (the fiber sum is denoted $\L' *_{\beta} \L''$) is 
constructed as follows: 

{\bf Step 1.} Take the disjoint union of metric graphs $(L', \ell')\u 
(L'', \ell'')$ 
and identify  $v$ and $\be(v)$ for all $v\in S'$. The result is the metric 
graph $(L, \ell)$.  

{\bf Step 2.} Let $W$ be the image in $L$ of $W'\u W''$, we let $W$ be the 
marking  of the resulting fiber sum $\L:= (L, \ell, W)$. 
If the images $Z', Z''$ of $W', W''$ are given, we define the vector $Z$ 
(the image of $W$) as the vector with the coordinates $\phi(w_j)$, where $w_j\in W$ 
and $\phi$ is in $C(\L', Z')$ or in $C(\L'', Z'')$. 

In the case $\L'= \L''$ the above construction has an analogue 
which we call  the {\em self-fiber sum}. The only difference is 
that on the first step instead of $L'\u L''$ we take the same 
graph $L'$ as before. The self-fiber sum will be denoted   $\L' *_{\beta}$. 

\begin{rem}
Notice that if $\L'\cong \L''$ then $\L' *_{\beta}$ is different from 
$\L' *_{\beta} \L''$.  
\end{rem}

If $S', S''$ are singletons $\{u\}, \{v\}$ then we will 
denote $\L'*_{\be}$ by $\L'*_{u=v}$. In what follows we will 
consider $\L', \L''$ being canonically embedded in $\L$. 

\subsection{Functional linkages}

Let $\k$ denote either $\C$ or $\R$. We will identify $\C$ with $\R^2$ 
and $\R$ with the real axis in $\C$. Recall that $C(\L,Z)$ is the set 
of real points of an affine scheme ${\mf C}(\L,Z)$. We now give the main definition. 

\begin{defn}
\label{2.7}
Let $\O\in \k^m$ and $F: \k^m \to \k^n$ be a map. We define a $\k$-functional 
linkage $\L$ for the germ $(F, \O)$ as follows: 

It is an abstract marked linkage 
$\L= (L, \ell, W)$ with $m$ distinguished vertices $P_1,.., P_m$ (called the 
{\bf input} vertices) and $n$ additional distinguished vertices $Q_1,..., Q_n$ 
(called the {\bf output} vertices) and a particular choice of a 
vector $Z\in \C^s$, the image of marking. We require this data to satisfy 
the axioms:

(1) The forgetful map $p:C(\L, Z)\to (\R^2)^m$ given by 
$$
p(\phi)= (\phi(P_1),..., \phi(P_m)), \quad \phi \in C(\L, Z)
$$
is a regular topological branched covering of a  
domain $Dom(\L, Z)$ in $\k^{m}$, so that the group $Sym(\L, Z)$ 
of automorphisms of $p$ consists of scheme-theoretic automorphisms. 
We let $Crit(\L, Z)$ denote the set of irregular points of the 
scheme-theoretic morphism $p$ and $C^*(\L, Z):= C(\L, Z) - Crit(\L, Z)$. 
It is clear that $C^*(\L, Z)$ is invariant under $Sym(\L, Z)$. 
We let $Dom^*(\L, Z):= p(C^*(\L, Z))$. Thus  
$$
p: C^*(\L, Z) \to Dom^*(\L, Z)
$$ 
is a locally analytically trivial\footnote{In the scheme-theoretic sense.}  
covering. We require $\O \in Dom^*(\L)$. 

(2) The forgetful map $q: C(\L, Z)\to \R^{2n}$ given by 
$$
q(\phi)= (\phi(Q_1),..., \phi(Q_n)), \quad \phi \in C(\L, Z)
$$
factors through $p$ and induces the map $F|Dom(\L,Z): Dom(\L,Z) \to \k^n$. 
We will say that the germ $(F, \O)$ is {\bf defined} by the linkage 
$\L$ and the vector $Z$.  
\end{defn}

Notice that in the definition of functional linkage for a germ $(F, \O)$ 
the metric ball around $\O$ which is contained in $Dom^*(\L, Z)$ is not specified. 
We will also need the following modification of the above definition:

\begin{defn}
Suppose that the pair $(\L, Z)$ as above defines the germ 
$(F, \O)$ and, moreover, $U$ is a \nbd of $\O$ such that 
$U\subset Dom^*(\L, Z)$. Then we say that the pair 
$(\L, Z)$ {\bf defines} $(F, U)$.  
\end{defn}

The group $Sym(\L,Z)$ will be called the {\em symmetry group} of $\L$ and 
$Dom(\L, Z)$ will be called the {\em domain} of $\L$ 
(of course they both depend on $Z$). The set of {\em input vertices} is denoted 
by $In(\L)$ and the set of output vertices by $Out(\L)$. 
We will refer to $\R$-functional linkages as {\em real functional 
linkages} and $\C$-functional linkages as {\em complex functional linkages}. 
If the choice of $\k, Z, In(\L), Out(\L)$, $\O, U$ or/and $F$ is suppressed 
then $\L$ is also referred as a functional linkage. 

\begin{lem}
Suppose that $\L$ is a functional linkage. 
A point $\phi\in C(\L,Z)$ belongs to $Crit(\L,Z)$ iff either 
$Dp_{\phi}: T_{\phi} {\mf C}(\L)\to T_{p(\phi)}Dom(\L,Z)$ has nonzero kernel or 
$\phi$ is not a smooth point of ${\mf C}(\L,Z)$. 
\end{lem}
\proof If $Dp_{\phi}$ has nonzero kernel then clearly $\phi\in Crit(\L,Z)$. 
If $\phi\notin Crit(\L,Z)$ then $\phi\notin Wall(\L)$ and hence $p(\phi)$ 
is in the interior of $Dom(\L,Z)$. If $p: ({\mf C}(\L,Z), \phi) \to (Dom(\L,Z), p(\phi))$ 
is an isomorphism of analytical germs then $({\mf C}(\L,Z), \phi)$ is 
necessarily smooth.   

Now suppose that $\phi\in Crit(\L,Z)$ and $Dp_{\phi}$ has zero kernel. 
(In particular $\phi\notin Wall(\L)$.) 
If $\phi$ is a smooth point then the dimension of $C(\L,Z)$ near $\phi$ 
is the dame as the dimension of the interior of $Dom(\L,Z)$, 
thus we can apply the implicit function theorem to conclude that 
$\phi\in C^*(\L,Z)$.  \qed  

\begin{rem}
Actually $\phi$ is irregular if and only if $Dp_{\phi}$ has nonzero kernel 
in the case when $\L$ is produced from {\bf elementary} linkages 
via composition (as in the 1-st and 2-nd 
functionality theorems), however we do not need this fact 
in the present paper and will not prove it here. 
\end{rem}

\begin{cor}
$Dom^*(\L,Z)$ and $C^*(\L,Z)$ are open subsets in $\k^{m}$ and $C(\L,Z)$ 
respectively. 
\end{cor}
\proof The set of singular points of ${\mf C}(\L,Z)$ is closed as well 
as the set of points where $Dp_{\phi}$ has nonzero kernel. The 
restriction of $p$ to $C^*(\L,Z)$ is open.  \qed   

\medskip
Besides functional linkages we will also need {\em closed functional linkages} 
defined as follows. Let $\L'$ be a functional linkage and $\be$ be a 
map from a subset of $Out(\L')$ to $W$ (the marking  
of $\L'$). Then the linkage $\L:= \L'*_{\be}$ is called a {\em closed functional 
linkage}. Such linkage still has the {\em input map} 
$p$ (the restriction of the input map $p'$ of $\L'$) and 
the {\em domain} $Dom(\L,Z)$ which is the image of $p$. 
The group of symmetries  $Sym(\L', Z)$ acts naturally on 
$C(\L, Z)$ and we let 
$$
Sym(\L, Z):= Sym(\L', Z)
$$

Let $\L$ be a (possibly closed) functional linkage. We have the  
scheme-theoretic input morphism $p$. As above we let 
$Crit(\L,Z)\subset C(\L,Z)$ denote the set of irregular points for $p$, 
it is invariant under the action of the symmetry group and we let  
$C^*(\L,Z):= C(\L,Z)- Crit(\L,Z)$, $Dom^*(\L,Z):= p(C^*(\L,Z))$. 
To check that $p$ is a local analytic 
isomorphism at $\phi$ it is enough to verify that $p$ induces 
bijections of Zariski tangent spaces of all orders: 
$T^{(n)}_{\phi}{\mf C}(\L,Z)\to T^{(n)}_{p(\phi)}{\mf Dom}(\L,Z)$, 
see Lemma \ref{algeo}.  

If $\L$ is a (possibly closed) functional linkage we let $Wall(\L,Z)\subset 
C(\L, Z)$ denote the collection of fixed points of nontrivial 
elements of $Sym(\L,Z)$. As we shall see (in \S \ref{funthm3}) $Wall(\L,Z)$ 
is nowhere dense in $C(\L, Z)$. 
Components of $Wall(\L,Z)$ will be called {\em walls} 
and components of $Crit(\L,Z)- Wall(\L,Z)$ will be called {\em quasiwalls}. 
We retain the names {\em walls} and {\em quasiwalls} for projections of 
walls and quasiwall to $Dom(\L,Z)$ via the input map $p$. 

\begin{rem}
\label{cri}
Suppose that $\L$ is a (possibly closed) functional linkage and 
$\phi\in C(\L, Z)$ is a relative 
realization such that all the derivatives 
$(Dp)^{N}_{\phi}: T_{\phi}^N{\mf C}(\L, Z)\to T_{p(\phi)}^N{\mf Dom}(\L, Z)$ 
are bijective. Then $\phi\in C^*(\L,Z)$. 
\end{rem}

\medskip
If $\L$ is a (possibly closed) functional linkage and $v$ is an input or 
output vertex of  $\L$ then we let $\pi_v: Dom(\L,Z)\to \k$ to denote the 
composition $eval_v\circ p^{-1}$. Clearly this composition is a 
well-defined $\k$-linear mapping.

\section{Functionality theorems}
\label{fun}

In this section we prove three theorems establishing that 
functionality of linkages is preserved (under appropriate circumstances) 
by the (self-) fiber sum. To simplify the notations in this section 
we will suppress the choices of fixed vertices $W$ and their images $Z$ 
under relative realizations.

Let $\k$ be either $\R$ of $\C$ and $\overline{\k}$ be its algebraic 
closure. Consider $\a, \b$ which are $\overline{\k}$- and 
$\k$-functional linkages respectively for the functions 
$f:\overline{\k}^n \to \overline{\k}^s$, 
$g:\k^m \to \k^t$. Let $T= \{P_1,..., P_t\}\subset In(\a)$ be  
 a collection of the input vertices and $Out(\b)= \{Q_1,..., Q_t\}$. 
We will assume 
that $t=n$ in the case $\k\ne \overline{\k}$. 
Suppose that we are given a bijection $\be: T\to Out(\b), 
\be(P_i)=Q_i$. Let $x=(x_1,..., x_n)$ and $y=(y_1,..., y_m)$ denote the 
coordinates in $\overline{\k}^n, \k^m$ respectively. 

The goal of the ``1-st functionality theorem''  below is to 
prove that the fiber sum $\L= \a *_{\be} \b$ is a functional linkage 
for the composition of the functions $f, g$ and to describe the 
configuration space, domain, etc.,  of the linkage $\L$. We  will refer to 
 $\L= \a *_{\be} \b$ as the {\em composition} of the linkages $\a, \b$. 

%\begin{rem}
%\label{dega}
%If we have a pair of input vertices in either $\a$ or $\b$ then 
%they are not connected by any edge. This implies that each degenerate 
%triangle in $\L$ is contained either in $\a$ or in $\b$. 
%\end{rem}

We let $eval_T, \pi_T$ be the vector-functions  
$$
eval_T:=(eval_{P_1},..., eval_{P_t}): C(\a) \to \k^t, \quad 
\pi_T:= (\pi_{P_1},..., \pi_{P_t})
$$
Here $n\ge t$ and we assume that $\k^t$ is canonically embedded 
in $\k^n=\k^t\times \k^{n-t}$. We let $p', p''$ be the input maps of 
$\a, \b$ and $q', q''$ be the output maps of $\a, \b$.

\begin{thm}
\label{funthm1}
(The 1-st functionality theorem.) Suppose that 
$$
int Dom(\b) \cap g^{-1} int Dom(\a)\ne \0
$$
Then:  
\begin{enumerate}
\item The scheme ${\mf C}(\L)$ is isomorphic to the fiber product 
${\mf C}(\a)\times_{eval_{T}= q''}{\mf C}(\b)$. 

\item $\L$ is a $\k$-functional linkage for the composition $h$ of 
the functions $f, g$:
$$
f(g_1(y),..., g_{t-1}(y), g_t(y), x_{t+1},..., x_n)
$$
and  $In(\L):= In(\a)\cup In(\b)- T$, $Out(\L):= Out(\a)$. $Dom(\L)$ 
is isomorphic $Dom(\a)\times_{\pi_{T}= g} Dom(\b)$ (as semi-algebraic sets).     
\item $Sym(\L)\cong Sym(\a)\times Sym(\b)$. 
\item $Wall(\L)= Wall(\a) \times_{eval_{T}= q''} C(\b) \cup 
C(\a) \times_{eval_{T}= q''} Wall(\b)$. 
\item $Crit(\L)\subset Crit(\a)\times_{eval_{T}= q''} C(\b) \cup 
C(\a) \times_{eval_{T}= q''} Crit(\b)$. 
\item If $O'\in Dom^*(\a), O''\in Dom^*(\b)$ then $(O', O'')\in Dom^*(\L)$. 
If $t=n$ then $Dom^*(\b) \cap g^{-1}(Dom^*(\a))\subset Dom^*(\L)$.  
\end{enumerate}
\end{thm}
\proof (1) The first assertion is proven similarly to \cite[Theorem 8.17]{KM6}.

\medskip
Now we start proving (2). 
We let $p, q$ be the input and output maps of $\L$ where 
$In(\L):= In(\a)\cup In(\b)- T$, $Out(\L):= Out(\a)$. 

We define the isomorphism 
$$
\theta: Dom(\a)\times_{\pi_T= g} Dom(\b)\to Dom(\L)= p(C(\L))
$$ 
of semi-algebraic sets by the formula:
$$
\theta: (x_1,..., x_n)\times (y_1,..., y_m)\mapsto 
(y_1,..., y_m, x_{t+1}, ..., x_n)$$
The inverse $\theta^{-1}$ is given by the formula:
$$
(y_1,..., y_m, x_{t+1}, ..., x_n)\mapsto 
(g(y), x_{t+1}, ..., x_n)\times (y_1,..., y_m) $$
The image of $p$ has nonempty interior by the assumption of theorem, 
i.e. it is a {\em domain}. We leave the proof of the  equality
$$
h= q\circ p^{-1}
$$
to the reader. 

To conclude the proof of (2) it remains to show that $p$ is a regular 
ramified covering and its group of automorphisms consists of polynomial 
automorphisms of $C(\L)$, this will follow from the proof of (3) and (4) 
below.

\medskip
(3) Notice that the group $Sym(\a)\times Sym(\b)$ acts naturally on
$C(\L)$. Indeed, if $\eta= (\eta', \eta'') \in C(\L)\subset C(\a)\times C(\b)$ 
and $\si= (\si',\si'')\in Sym(\a)\times Sym(\b)$ then 
$$
\si(\eta):= (\si'(\eta'), \si''(\eta''))\in C(\a)\times C(\b)
$$
However $\si'(\eta')|T= \eta'|T$,  $\si''(\eta'')|In(\b)= \eta''|In(\b)$ since 
$\si', \si''\in Sym(\a), Sym(\b)$. Hence 
$$
(\si'(\eta'), \si''(\eta''))\in C(\L)= C(\a)\times_{eval_{T}= q'}C(\b)
\subset C(\a)\times C(\b)
$$
This implies that $\si$ acts on $C(\L)$. We leave it to the reader to 
verify that the action is faithful. It is clear that for each 
$\si\in Sym(\a)\times Sym(\b)$ we have:
$$
p\circ \si= p, \quad q\circ \si= q
$$
Suppose that $\phi,\psi\in C(\L)$ are such that 
$\phi|In(\L)= \psi|In(\L)$. Then 
$$
\phi|In(\b)= \psi|In(\b)
$$
which implies that $\phi|T= \psi|T$ since $\b$ is functional. Therefore 
$$
\phi|In(\a)= \psi|In(\a)
$$
It follows that there are symmetries $\si'\in Sym(\a), \si''\in Sym(\b)$ 
such that
$$
\si'(\phi|\a )=  \psi|\a , \quad \si''(\phi|\b )= \psi|\b  
$$
We conclude that $\psi= (\si' \times\si'')(\phi)$, in particular
$$
\psi|Out(\L)= \phi|Out(\L)
$$
This also shows that $Sym(\L)\subset Sym(\a)\times Sym(\b)$. 
This finishes the proof of (3). 

(4) Suppose that $\si= \si'\times\si''\in Sym(\L)$. 
Then $Fix(\si)$ is contained in the union
$$
Fix(\si'\times 1)\cup Fix(1\times\si'')
$$
Note that 
$$
Fix(\si'\times 1)= Fix(\si')\times_{eval_T= q'} C(\b)$$
$$
 Fix(1\times\si'')= C(\a) \times_{eval_T= q'} Fix(\si'')
$$
This proves (4). We leave it to the reader to verify that 
$p|C(\L)-Wall(\L)$ is a local homeomorphism into $Dom(\L)$. Since 
$p(\phi)=\psi$ iff there is $\si\in Sym(\L)$ such that $\psi= \si(\phi)$ 
we conclude that $p|C(\L)-Wall(\L)$ is a covering onto its image. 
This concludes the proof of (2).

(5) A realization (of a functional linkage) is an irregular point 
of the input map iff some of the $N$-th order derivatives of the 
input map at this point is not a bijection (see Remark \ref{cri}). 
Suppose that $\phi= (\phi',\phi'')\in C(\L)$ is such that 
$\phi'\in C^*(\a), \phi''\in C^*(\b)$  (i.e. both are not 
irregular), $x'=p'(\phi')$, $x''=p''(\phi'')$.  

Suppose that there are  $N$-th order infinitesimal deformations 
$\xi, \eta$ of $\phi$ such that:
$$
\xi|In(\b)= \eta|In(\b), \quad \xi|(In(\a)- T) = \eta|(In(\a)- T)
$$ 
If $\xi|\b\ne \eta|\b$ then $\phi''\in Crit(\b)$ which contradicts 
our assumption. If $\xi|\b=\eta|\b$ then 
$\xi|T=\eta|T$ and $\xi|In(\a)=\eta|In(\a)$ which implies that 
$\phi'\in Crit(\a)$. This contradiction shows that $D^N_{\phi}p$ is injective.  

To show surjectivity pick two $N$-th order tangent vectors 
$\tau'\in T^{N}_{x'}Dom(\a), \tau''\in T^{N}_{x''}Dom(\b)$ such that
$$
D^{N}\pi_T(\tau')= D^Ng(\tau'')
$$
Since $\phi''\notin Crit(\b)$, we can lift $\tau''$ to an $N$-th order 
tangent vector $\xi''\in T^N_{\phi''} {\mf C}(\b)$. Necessarily 
$\xi''|Out(\b)=\tau'|T$. By the same reason, $\tau'$ lifts to 
$\xi'\in T^N_{\phi'} {\mf C}(\a)$ and $\xi'|T = \tau'|T= \xi''|Out(\b)$. Thus 
$(\xi',\xi'')\in T^{N}_{\phi}{\mf C}(\L)$ since the latter is the fiber product 
of $T^N_{\phi'}{\mf C}(\a)$ and $T^N_{\phi''}{\mf C}(\b)$. 

This proves that 
$$
Crit(\L)\subset Crit(\a)\times_{eval_{T}= q'} C(\b) \cup 
C(\a) \times_{eval_{T}= q'} Crit(\b)
$$

(6) Follows directly from (5).  $\qed$ 

\medskip
Now we consider the self-fiber sums: $\L= \a*_{\be}$, 
where the linkage $\a$ is $\k$-functional for a vector-function 
$f(x_1,..., x_n)$ with the components $(f_1,..., f_m)$.  

\medskip
{\bf Case I:} $\be: S'=\{v\}\subset In(\a) \to In(\a)$, $\be(v)=w$.  
Then 
$${\mf C}(\L)\cong \{\phi\in {\mf C}(\a)|\phi(v)=\phi(w)\}$$
where $\cong$ is a scheme-theoretic isomorphism. 
Let $p'$ denote the input map of $\a$. Consider the input mapping
$$
p: C(\L)\to Dom(\L):= p(C(\L))
$$
which is the restriction of $p'$. Then $Dom(\L)$ equals
$$
\{ x\in Dom(\a) : \pi_v(x)= \pi_w(x)\}
$$
i.e. the intersection of $Dom(\a)$ with a hyperplane. We will assume that 
the intersection of $Dom^*(\a)$ with this hyperplane is nonempty.  Since  
$Dom^*(\a)$ is open we conclude that $Dom(\L)$ has nonempty interior. 

To identify 
the collection of critical points $Crit(\L)$ of $\L$ let $\phi\in C(\L)$ 
be the restriction of $\phi'\in C(\a)$. Suppose that $\xi, \eta$ are distinct 
$N$-th order infinitesimal deformations of $\phi$ 
such that $\xi|In(\L)= \eta|In(\L)$. Then $\xi,\eta$ lift to $N$-th order 
infinitesimal deformations $\xi',\eta'$ of $\phi'$ such that 
$\xi'|In(\a)=\eta'|In(\a)$. This implies that $\phi'\in Crit(\a)$. 

Suppose now that $Dp_{\phi}^N$ is not surjective. Then we can find 
$$
\xi\in T^N_{p(\phi)}Dom(\L)
$$
which can not be lifted to $T^N_{\phi} {\mf C}(\L)$. If $\xi$ can not be lifted 
to  $T^N_{\phi} {\mf C}(\a)$ then $\phi'\in Crit(\a)$ and we are done. 
Suppose that we can find a lift $\t\xi\in T^N_{\phi} {\mf C}(\a)$ of $\xi$. 
Note that $\t\xi|v =\xi|v= \xi|w=\t\xi|w$. Thus 
$\t\xi\in T^N_{\phi} {\mf C}(\L)$. 
We  conclude that $Crit(\L)\supset Crit(\a)\cap C(\L)$. 

The group of symmetries $Sym(\a)$ acts naturally on $C(\L)$. Let 
$\phi, \psi\in C(\L)$ be realizations such that 
$\phi|In(\L)= \psi|In(\L)$. Then $\phi, \psi$ lift to realizations 
$\phi', \psi'\in C(\a)$ such that $\phi'|In(\a)= \psi'|In(\a)$. It follows that 
there is a symmetry $\si'\in Sym(\a)$ such that $\si'\psi'= \phi'$, hence 
$\si\psi= \phi$, where $\si$ is the restriction of $\si'$ to $C(\L)$. 
We conclude that $Sym(\L)\cong Sym(\a)$. Thus we have proved 

\begin{thm}
\label{funthm2}
(The 2-nd functionality theorem.) Suppose that the set 
$\{ x\in Dom^*(\a) : \pi_v(x)= \pi_w(x)\}$ 
is nonempty. Then: the mapping $p: C(\L)\to Dom(\L)$ is a regular 
ramified covering with the group $Sym(\L)\cong Sym(\a)$ of covering 
transformations. This covering is a (scheme-theoretic) locally 
analytically trivial covering over 
$$
Dom^*(\L)\supset Dom(\L) \cap Dom^*(\a).$$ 
The real semi-algebraic set $Dom(\L)$ is isomorphic to 
$\{x\in Dom(\a): \pi_v(x)=\pi_w(x)\}$.  The linkage $\L$ is 
functional for the restriction of the function $f$ to the 
hyperplane %\newline 
$\{\pi_v(x)=\pi_w(x)\}$. 
\end{thm}

{\bf Case II:} $\be: S'\subset Out(\a)\to W$, whence 
the linkage $\L= \a*_{\be}$ is a ``closed functional linkage''. 
Let $S'= \{Q_1,..., Q_t\}$, then $\be(Q_j)$ is a fixed vertex 
for each $j=1,...,t$. Let $z_j:= \phi(\be(P_j))$ for all relative 
realizations. 
Consider the input map for $\L$:
$$
p: C(\L)\to Dom(\a)
$$
which is the restriction of the input map $p'$ of $\a$. 
The image of this restriction is
$$
Dom(\L)= \{ x\in Dom(\a) : f_j(x)=z_j, j=1,..., t\}
$$
Now we describe the set of critical points 
of the mapping $p$. Suppose that $\phi\in C(\L)$ is a realization which 
corresponds to $\phi'\in C^*(\a)$ and $\xi, \zeta$ are infinitesimal 
deformation of $\phi$ (of the order $N$) such that
$$\xi|In(\L)=\zeta| In(\L)$$
Then $\xi, \zeta$ determine $N$-th order infinitesimal deformations $\xi', \zeta'$ 
of $\phi'$. Clearly $\xi'|In(\a)=\zeta'|In(\a)$. This implies that 
$\xi'=\zeta'$, which means that for each $\phi\in C^*(\L)$ the 
push-forward map
$$
D^{(N)}p: T^{(N)}_{\phi} {\mf C}(\L) \to T^{(N)}_{p(\phi)}{\mf Dom}(\L)
$$
is injective. To show surjectivity we pick a vector $\eta\in  
T^{(N)}_{p(\phi)}{\mf Dom}(\L)$.  Since $\phi\notin Crit(\a)$, we can lift $\eta$ to  
$\t\eta\in T^{(N)}_{\phi} {\mf C}(\a)$. The vector $\t\eta$ belongs to 
 $T^{(N)}_{\phi} {\mf C}(\L)$ iff $\t\eta|Q_i=0$ for all $i=1,..., t$. However
$(\t\eta|Q_1,..., \t\eta|Q_t)$ equals $(D^{(N)}f)\eta$ since $\a$ is functional 
and $\phi\notin Crit(\a)$.  
Since $\eta$ is tangent to $Dom(\L)$ we conclude that  
$D^{(N)}f(\eta)=0$, which implies that 
$$
\t\eta \in T^{(N)}_{\phi} {\mf C}(\L)
$$
Surjectivity of $D^{(N)}p$ at $\phi$ follows. We conclude that 
$$
p: C^*(\L)\to Dom(\L)
$$ 
is a local analytic isomorphism onto its image. Now consider the 
group of symmetries of $\L$. Suppose that $\phi, \psi\in C(\L)$ 
are such that $\phi|In(\L)= \psi|In(\L)$. It follows that they are 
restrictions of realizations $\phi', \psi'\in C(\a)$ and there is a symmetry 
$\si\in Sym(\a)$ such that $\phi'= \si\psi'$. 
Therefore the mapping $p$ is an analytically trivial covering over 
$Dom^*(\L)$ with the group $Sym(\L)\cong Sym(\a)$ of covering 
transformations. We have proved 

\begin{thm}
\label{funthm3}
(The 3-rd functionality theorem.) 
Suppose that $Dom(\L)\cap Dom^*(\a)$ is nonempty. 
Then the input mapping $p: C(\L)\to Dom(\L)$ is a regular ramified 
covering with the group $Sym(\L)\cong Sym(\a)$ of covering 
transformations. This covering is locally analytically trivial 
over $Dom^*(\L):= Dom(\L) \cap Dom^*(\a)$ (in the scheme-theoretic sense). 
The set $Dom^*(\L)$ (which is an open\footnote{In the classical topology.} 
subset in $Dom(\L)$) is analytically isomorphic to 
$$\{x\in Dom^*(\a): f_j(x)=z_j, j=1,...,t\}\quad .$$  
\end{thm}

\section{Fixing fixed vertices}
\label{fixing}

The goal of this section is to relate relative configuration spaces $C(\L, Z)$ 
of marked linkages and the moduli spaces $\M(\L)$ of based linkages. 
Let $\L= (L, \ell, W)$ be a marked linkage, 
$Z= (z_1,..., z_s)\in \C^s$ and $W= (w_1,..., w_s)$. Pick any relative 
realization $\phi\in C(\L, Z)$. 

We first let $\L'$ be the disjoint union of $\L$ and the metric 
graph $\I$ which consists of a single edge $e^*$ of the unit length 
connecting the vertices $v_1, v_2$. Choose the isometric embedding 
$\phi= \phi_{\I}: \I\to \C$ which maps 
$v_1$ to $0$ and $v_2$ to $1\in \R$. We get a map $\phi: W\cup \V(\I) \to \C$. 
Then for each pair of vertices 
$a, b\in W\cup \V(\I)$ we do the following: 

(a) If $\phi(a)= \phi(b)$ for $\phi\in C(\L, Z)$, we identify the vertices $a, b$.   

(b) Otherwise add to $\L'$ the edge $[a b]$ of the length 
$|\phi(a)- \phi(b)|$. 

\medskip
\noindent Let $\t\L$ be the resulting based linkage (with the distinguished 
edge $e^*=[v_1 v_2]\subset \I$). 

There are now two different cases: (i) the vector $Z$ is not real (we shall assume 
$z_1\notin \R$), (ii) $Z$ is real. 
In the second case for each realization $\phi\in \M(\t\L)$ we have: 
$$
(\phi(w_1), \phi(w_2),..., \phi(w_s))= (z_1,..., z_s)\in \R^s
$$
Thus the natural (scheme-theoretic) morphism $\iota: C(\L, Z)\to \M(\t\L)$ 
is a bijection. 
However (unless the image of $W\cup \V(\I)$ in $\t\L$ consists of 
two vertices) we had created new nilpotent elements in the ring of $\M(\t\L)$, thus 
$\iota$ is not a scheme-theoretic analytic isomorphism. On 
the other hand, since we are interested in the real 
reduced schemes, we can use the 
same trick as in the case of rigidified parallelograms: we give 
$\M(\t\L)$ the scheme-theoretic structure of ${\mf C}(\L, Z)$. 

In the case (i)  for each realization $\phi\in \M(\t\L)$ we have: 
$$
(\phi(w_1), \phi(w_2 ),..., \phi(w_s))= (z_1,..., z_s) 
$$
or
$$
(\phi(w_1), \phi(w_2),..., \phi(w_s))= (\bar{z}_1,..., \bar{z}_s) \ne (z_1,..., z_s)
$$
On the other hand, we did not create new nilpotent element in the ring of 
$\M(\t\L)$
(since for each $w_i$ (if $i\ge 2, z_i\ne 0$) either the triangle 
$\Del(v_1 v_2 w_i)$ or 
the triangle $\Del(v_1 w_1 w_i)$ is nondegenerate. Thus in the case (i) 
we get an analytically trivial covering  $\tau: \M(\t\L)\to C(\L, Z)$ 
given by:  
 
For $\phi \in \M(\t\L)$ we let $\tau(\phi):= \phi|\L$ if 
$\phi(w_1)=  z_1$ and $\tau(\phi):= \overline{\phi}|\L$ if 
$\phi(w_1)=  \overline{z_1}$. 

\no This covering has a section $\si:\psi\in C(\L, Z)\mapsto \M(\t\L)$ 
such that: 
$$
\si(\psi)|\L := \psi, \quad \si(\psi)|\I:= \phi_{\I}
$$
It is clear that the group of automorphisms of the covering $\tau$ is 
$\Z_2$ and is generated by the complex conjugation. We summarize this in 
the following

\begin{lem}
\label{fi}
(i) In the case $Z\notin \R^s$ there is an 2-fold analytically 
trivial covering $\tau: \M(\t\L)\to C(\L, Z)$. 

(ii) In the case $Z\in \R^s$ there is an  
isomorphism $\tau: \M(\t\L)\to C(\L, Z)$.
\end{lem}

\section{Elementary linkages}
\label{elementary}

In this section we construct several {\em elementary} functional 
linkages: translators (for  the translation $z\mapsto z+b$), the adder 
(for the summation $(z,w)\mapsto z+w$), pantographs (for the functions  
$z\mapsto \la z$ and $z\mapsto -z$), inversors (for the functions 
$z\mapsto t^2/\bar{z}$), the multiplier\footnote{Strictly speaking, our 
multiplier is not so elementary, for instance it would be difficult to 
draw a picture of the corresponding graph.} (for the function $(z,w)\mapsto zw$) 
and the linkage for straight line motion. These linkages serve as building 
blocks for the proof of Theorem A. We make the following convention concerning 
usage of {\em elementary linkages}:

All elementary linkages come with parameters which do not affect functions 
that they define but affect domains of the linkages. Thus if we use 
several elementary linkages with the same name ${\cal N}$ in 
constructing another linkage $\L$ via fiber sum, we allow different 
choices of the 
parameters for different appearances of ${\cal N}$ in this fiber sum. 

We also omit the image of marking $Z$ in the notation for $Dom, Dom^*$ 
of the elementary linkages. 

All {\em elementary} linkages in the section (with the exception of the 
multiplier) are modifications of classical constructions, where 
appropriate modification was made 
to ensure functionality. We decided to avoid Kempe's construction of the 
multiplier \cite{Kempe} since computation of $Dom$ and $Dom^*$ for Kempe's 
linkage presents some difficulties, we use an algebraic trick instead.

\subsection{The translators}

\begin{figure}[tbh]
\leavevmode
\centerline{\epsfxsize=3in\epsfbox{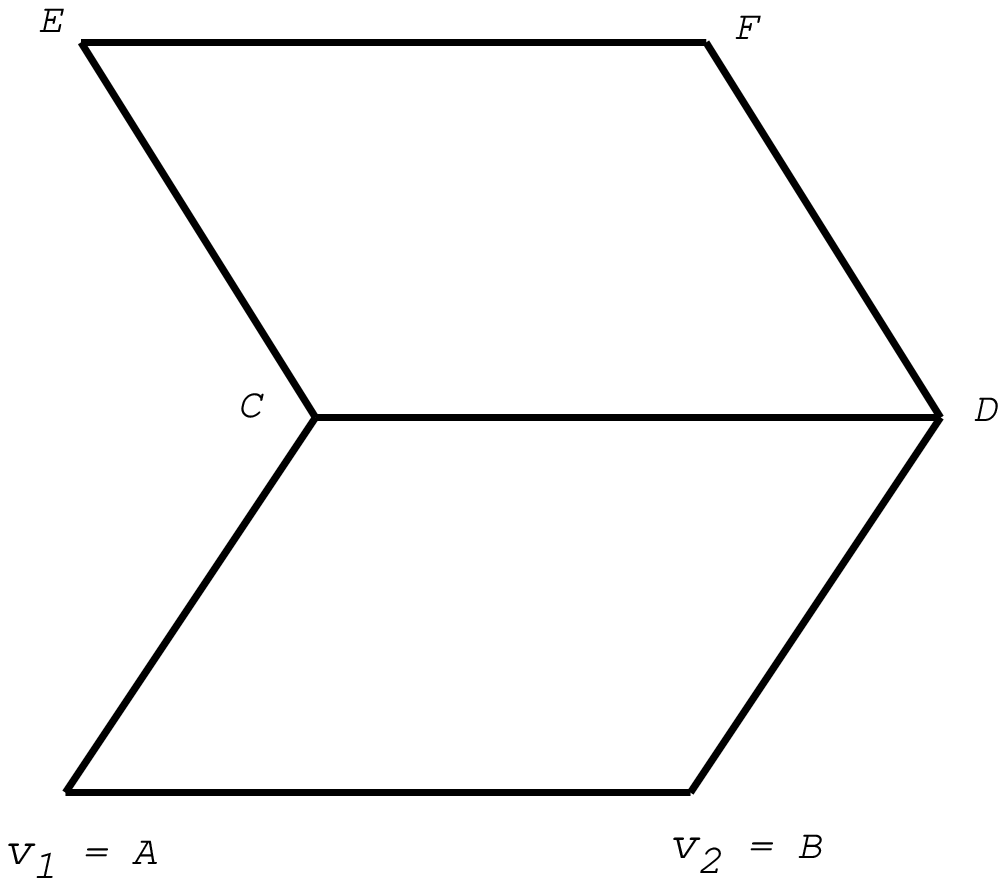}}
\caption{\sl A translator. The parallelograms $[ACDB]$ and $[CEFD]$ are 
rigidified. The set of intput/output vertices is $\{E, F \}$ and 
$s= \ell[AC]> t=\ell[CE]$.} 
\label{Fig10}
\end{figure}

Let $b$ be a fixed nonzero complex number. 
The  translation operations $\tau_b: z\mapsto z+ b$, $\tau_{b}':w\mapsto w-b$  
are defined using the translator which is described on Figure \ref{Fig10}. 
Depending on the operation either $F$ or $E$ is the input (resp. output). 
The point is that if $E$ is the input then by adjusting side-lengths of the 
corresponding translator $\T_{b}$ we can get any $z\in \C-\{0\}$ into 
$Dom^*(\T_{b})$. To get $0\in \C$ 
into $Dom^*$ we use the point $F$ as the input (and $E$ as the output) 
of a functional linkage $\T'_{b}$ for $\tau_{b}'$. Below we 
present the details. First of all let $W:= (v_1, v_2)$ be the marking 
of $\L$ (which is either $\T_b$ or 
$\T'_b$) and $Z:= (0, b)\in \C^2$. Below we shall use the 
relative configuration spaces of $\L$ associated to this data. 

The next lemma follows from the triangle inequalities and its proof is 
left to the reader. 

\begin{lem}
$Dom(\T_b)$ and $Dom(\T'_b)$ are the annuli 
given by the inequalities:
$$
Dom(\T_b)= \{\rho:=  s- t\le |\phi(E)|\le R:= s+t\}
$$
$$
Dom(\T'_b)=\{\rho\le |\phi(F)- b|\le R\}\ \qed 
$$
\end{lem}
Notice that the centers of these annuli are at the points $0, b$ 
respectively. 

\begin{lem}
Both linkages $\T_b$, $\T'_b$ are functional for 
the functions $\tau_b, \tau'_b$. Walls in the domains 
of these linkages are boundary circles of the corresponding annuli. 
\end{lem}
\proof We consider the first linkage $\L= \T_b$, 
the proof for the second linkage is analogous. Notice that 
for each realization $\phi$ we have: 
$$
\phi(F)= \phi(E)+ \phi(B), \quad \phi(D)= \phi(C)+ \phi(B)
$$ 
and $\phi(A)\ne \phi(E)$ (see Lemma \ref{2.6}). 
Suppose that $\phi\ne\psi$ are relative realizations and $\phi(E)= \psi(E)$. 
Then $\phi(F)= \psi(F)$. Thus $\psi(C)$ is the reflection of $\phi(C)$ 
in the line through   $\phi(A), \phi(E)$. Therefore $\L$  is functional 
and the group of symmetries of $\L$ is $\Z_2$. The fixed points of the 
generator of $Sym(\L)$ correspond to realizations for which the triangle 
$$
\Del(\phi(A),\phi(C),\phi(E))$$
 is degenerate, i.e. one of the inequalities defining $Dom(\L)$ is the  
equality. $\qed$

\medskip
Suppose that $\phi\in C(\L,Z)-Wall(\L,Z)$ 
is a realization of $\L= \T_b$ or $\L=T_b'$ 
such that none of the parallelograms $[\phi(A)\phi(B)\phi(D)\phi(C)]$, $[\phi(C)\phi(D)\phi(F)\phi(E)]$ is degenerate.  
We leave it to the reader to verify (using Lemma \ref{tri}) that in 
this case the morphism  
of analytic germs $p:({\mf C}(\L,Z), \phi)\to ({\mf Dom}(\L,Z), p(\phi))$ 
is invertible and hence such $\phi$ belongs to $C^*(\L,Z)$. For instance, 
if $F$ is the input then one can recover $\phi(D)$ are analytic function of 
$\phi(F)$ and $\phi(B)=b$ (since the triangle $\Del(\phi(B)\phi(D)\phi(F))$ 
is nondegenerate). Then one recovers $\phi(C)$ as analytic function of 
$\phi(A)=0$ and $\phi(D)$ (since the first of the parallelograms is nondegenerate), etc. 

On the other hand, if one of the above parallelograms is degenerate then 
the derivative
$$
Dp_{\phi}: T_{\phi}{\mf C}(\L,Z)\to {\mf Dom}(\L,Z)
$$
has zero kernel.  

We summarize this in the following lemma

\begin{figure}[tbh]
\leavevmode
\centerline{\epsfxsize=3in\epsfbox{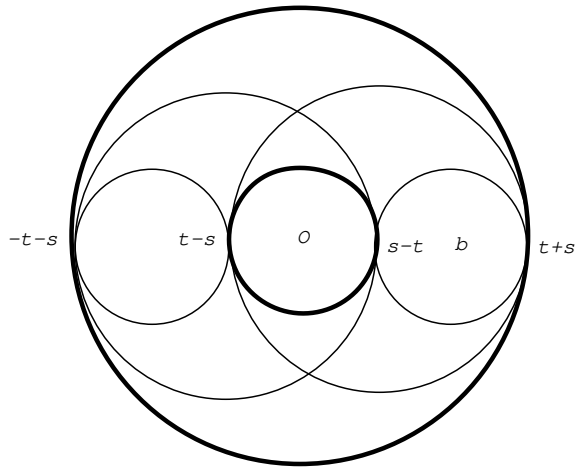}}
\caption{\sl Quasiwalls in the domain of the translator $\T_b$.} 
\label{f20}
\end{figure}

\begin{lem}
If $\phi\in C(\L,Z)$ belongs to a quasiwall then one of the parallelogram 
$$
[\phi(A) \phi(B) \phi(D) \phi(C)], \quad [\phi(E) \phi(F) \phi(D) \phi(C)]
$$
is degenerate. Let $\be:= b/|b|$. 
The closure of the union of quasiwalls for both $\T_b$, $\T_b'$ 
is the union of four circles:
$$
\{\phi(C)= \pm s \be, |\phi(E)\mp s\be|= t\}\cup
\{|\phi(E)\pm t\be|= s\}
$$
Unions of quasiwalls for both $\L$ are disjoint from the collection of 
$\phi$ for which $\phi(E)$ belongs to the line through $0, b$. 
\end{lem}

In particular, if 
$$s+t> |b| > s-t >0, \quad b\in \R$$
 then the origin belongs to $Dom^*(\T'_b)$.

\begin{rem}
In the following sections (except in \S \ref{draw}) each time when we 
have to use a translator we shall pick real numbers $b$. 
\end{rem}

The quasiwalls for $\T_b$ are described on Figure \ref{f20}, the 
picture for quasiwalls of $\T'_b$ is obtained via translating by $b$. 
Notice that
$$
int(B_t(\pm s))\subset Dom^*(\T_b), \quad 
int(B_t(\pm s +b)) \subset Dom^*(\T'_b)
$$
(provided that $b\in \R$). 

\begin{cor}
\label{trans}
Let $0< r < |b|, b\in \R$. Then the parameters  $t, s$ can be chosen so that
$$
B_r(b)\subset Dom^*(\T_b),  \quad 
B_r(0)\subset Dom^*(\T'_b)
$$
\end{cor}
\proof Choose $t, s$ so that $|t-s|\to 0$, $|t+s|\to \8$. Then 
$Dom^*(\T_b)$ converges to the union of half-planes $\{Re(z)\ne 0\}$ and  
 $Dom^*(\T'_b)$ converges to $\{Re(z)\ne b\}$.  \qed  

\subsection{The pantograph}
\label{pantograph}

\begin{figure}[tbh]
\leavevmode
\centerline{\epsfxsize=3in\epsfbox{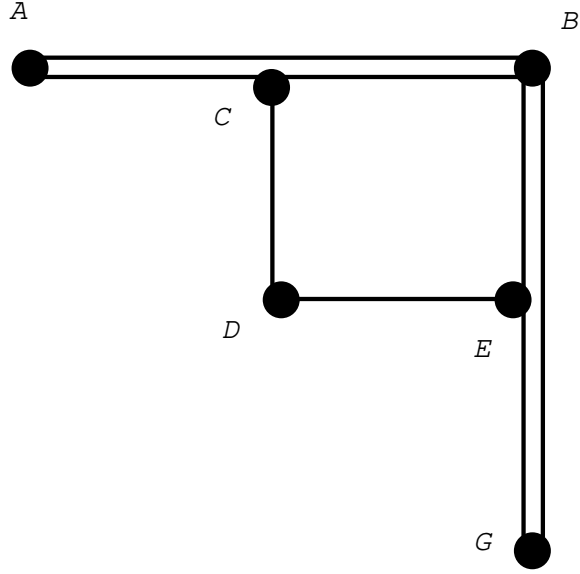}}
\caption{\sl The rigidified pantograph $\p$: 
the parallelogram $BCDE$ is rigidified, $\la >1$. This 
linkage is not marked, we shall use different choices of 
input/output vertices later on. We take: 
$s= \ell[AB]= \la \ell[AC] \ne t= \ell[BG]= \la \ell[BE]$.}
\label{Fig9}
\end{figure}

The (rigidified) pantograph $\p$ is described on Figure \ref{Fig9}, 
recall that we use the Convention \ref{contri} for the two degenerate 
triangles in $\p$ as well as for the rigidified parallelogram. 
The pantograph is a versatile linkage, its role in 
engineering\footnote{That goes 
back to at least 17-th century, see \cite{Sch}.}  was as a 
functional linkage for the functions $z\mapsto \la z, z\mapsto \la^{-1} z$, 
$\la >1$. 

\begin{rem}
In the next  section we shall also use the pantograph to construct the adder. 
\end{rem}

In the case of the function $z\mapsto \la z$ we let 
$W:= \{A\}$ be the fixed vertex, $Z:= 0$, take $D$ 
as input and $G$ as output, let $\p_{\la}$ be the resulting linkage 
(it will be functional for $z\mapsto \la z$). 
By switching input and output we obtain a 
functional linkage $\p_{1/\la}$ for $z\mapsto  z/\la$. 

By letting $\{D\}= W$ instead of $A$, the same $Z$ as before,  $\la=2$ 
and taking $A$ as input and $G$ as output we obtain a functional linkage 
for the function $z\mapsto -z$ in the complex plane. Notice that the 
condition $s\ne t$ implies that for each realization $\phi$  the points 
$\phi(A)$, $\phi(D)$, $\phi(G)$ are pairwise distinct. 

Below we describe $Dom$ and $Dom^*$ of the pantograph, the proofs are 
similar to the previous section and are left to the reader. 

\begin{lem}
\label{pant}
For each choice of the fixed vertex and input/output vertices described above,  
the pantograph is a functional linkage. The group of symmetries is 
isomorphic to $\Z/2$ and is generated by the reflection of $\phi(B)$ 
in the line $\nu$ through the points $\phi(A)$, $\phi(D)$, $\phi(G)$. 
The walls are described by the condition: $\phi(B)\in \nu$. There are no 
quasiwalls. $Dom^*(\p)$ is the interior of $Dom(\p)$. 
If $D$ or $A$ is the input and $G$ is the output 
then the domain of the pantograph is the annulus given by the inequalities  
$$
|s-t|/\la \le |\phi(D)-\phi(A)|\le (s+t)/\la
$$
If $G$ is the input and $D$ is the output  
then the domain of the pantograph is the annulus given by the inequalities 
$$
|s-t|\le |\phi(G)-\phi(A)|\le s+t
$$
For given $\la$ as  $|s-t|\to 0$  and $|s+t|\to \8$ the domains 
$Dom^*(\p)$ are convergent to punctured complex planes.  
\end{lem}

Note that zero does not belong to $Dom^*$ of the pantograph for any 
choice of input/output points. To resolve this problem we compose 
$\p$ with the appropriate translators:
$$
-z = -(z+b) + b= \tau_{b}( -\tau_{-b}'(z)) 
$$ 
$$
\la z = \la (z+b) - \la b= \tau_{-\la b}( \la \tau_{-b}'(z))
$$
$$
z/\la = (z+b)/\la - b/\la = \tau_{-b/\la }( \tau_{-b}'(z) /\la )
$$
where $b\in \R-\{0\}$. We call the linkages computing the above functions the 
{\em modified pantographs} and denote them 
$\p_{-}', \p_{\la}', \p_{1/\la}'$ respectively. We would like 
$Dom^*(\p_{\la}'), Dom^*(\p_{1/\la}')$ to contain arbitrarily large 
compacts. This is done as follows: 

\begin{lem}
Fix $\la > 1$ and let $r> 0$. 
Then we can choose $b\in \R$ and edge-lengths for the translators and for 
the pantographs $\p_{\la}, \p_{1/\la}$ so that
$$
B_r(0)\subset Dom^*(\p'_{\la}), \quad B_r(0)\subset Dom^*(\p'_{1/\la})
$$
\end{lem}
\proof We consider the case of the modified pantograph $\L:= \p'_{\la}$, 
the second case is similar. By Theorem  \ref{funthm1} we have: 
$$
Dom^*(\L) \supset Dom^*(\T_{-b}')\cap [Dom^*(\p_{\la})- b] \cap 
[\la^{-1} Dom^*(\T_{-\la b}) -b] 
$$
Below we analyse the triple intersection:

\no (1) By choosing appropriately the parameters $t, b$ in $\T_{-b}'$ we can 
guarantee that $Dom^*(\T_{-b}')$ contains arbitrarily large discs 
around the origin (see Corollary \ref{trans}). In particular, 
$B_r(0)\subset Dom^*(\T_{-b}')$. 

(2) The domain $[Dom^*(\p_{\la})- b]$ is obtained by translating 
$Dom^*(\p_{\la})$ by $-b$. Recall that $Dom^*(\p_{\la})$ is an open annulus 
centered at zero. By adjusting parameters in $\p_{\la}$ 
(and keeping $\la$ fixed) we can guarantee that $Dom^*(\p_{\la})$ 
contains the disk $B_{\rho}(b)$ for each $\rho< |b|$, see Lemma \ref{pant}. 
We conclude that if $|b|> r$ (and under appropriate choice of edge-lengths 
in $\p_{\la}$,  $\T_{-b}'$) the domain
$$
 Dom^*(\T_{-b}')\cap [Dom^*(\p_{\la})- b]
$$
contains the disk $B_r(0)$. 

(3) Lastly we consider the domain $\la^{-1} Dom^*(\T_{-\la b}) -b$, 
it contains the disk $B_r(0), r= R/\la,$ provided that 
$Dom^*(\T_{-\la b})$ contains the disk of radius $R$ centered at 
the point $\la b$. The domain $Dom^*(\T_{-\la b})$ is again an annulus 
centered at zero. By adjusting parameters of the linkage $\T_{-\la b}$ 
(and keeping $\la, b$ fixed) we can guarantee that $Dom^*(\T_{-\la b})$ 
contains arbitrary annuli centered at zero. Hence for each $R< |\la b|$ 
(i.e. $r < |b|$) and under appropriate choice of edge-lengths in 
$\T_{-\la b}$, the domain $Dom^*(\T_{-\la b})$ contains 
the disk $B_R(\la b)$. We conclude that for each $\la > 0$ and $r> 0$ 
if we choose $b$ such that $r < |b|$, then the edge-lengths of the linkage $\L$ 
can be chosen so that 
$$
B_r(0) \subset Dom^*(\L).\quad\quad\quad \qed
$$

\subsection{The adder}

We again consider the rigidified pantograph, only now $\la=2$, 
the vertices $A, G$ are the inputs, $D$ is the output and there 
is no fixed vertices at all. We will use the notation ${\mathcal Q}$ 
for the resulting linkage. Similarly to the previous section ${\mathcal Q}$ 
is $\C$-functional for the function
$$
(z,w)\mapsto (z+w)/2
$$
(the input $A$ corresponds to $z$ and the input $G$ corresponds to $w$). 
As before the domain of ${\mathcal Q}$ is given by 
$$
\{(z,w)\in \C^2: t\le |z-w| \le 3t\}
$$
and $Dom^*({\mathcal Q})$ is the interior of $Dom({\mathcal Q})$. Note  
that the point $(b,-b)$ belongs to $Dom^*({\mathcal Q})$ provided that
$t< 2|b| < 3t$. The point $(0,0)$ does not belong to $Dom^*({\mathcal Q})$, 
similarly to the previous section we use appropriate translators to resolve 
this problem:
$$
(z+w)/2 = [(z+b) + (z-b)]/2 
$$ 
where $t< 2|b| < 3t$. Thus we get a modified linkage ${\mathcal Q}'$ for 
$(z,w)\mapsto (z+w)/2$ such that $(0,0)\in Dom^*({\mathcal Q}')$. 
To get the linkage $\L_A$ for the addition we combine ${\mathcal Q}'$  
and the modified pantograph $\p'_{2}$ for the multiplication by $2$: 
$$
(z,w) \mapsto (z+w)/2 \mapsto z+w
$$ 
Then $(0,0)\in Dom^*(\L_A)$.

\subsection{The modified inversor}

\begin{figure}[tbh]
\leavevmode
\centerline{\epsfxsize=3in\epsfbox{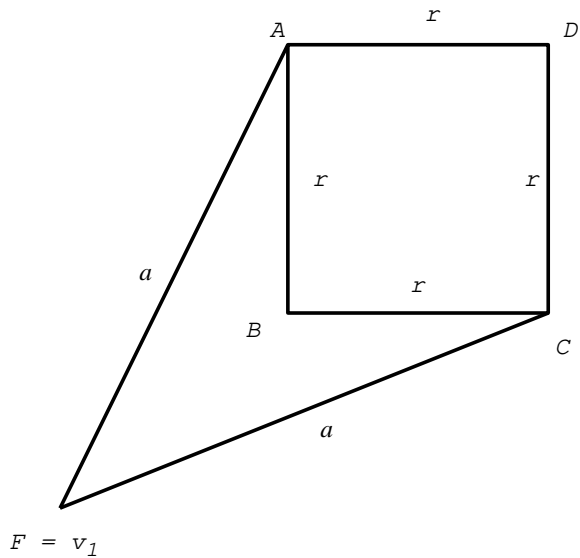}}
\caption{\sl The Peaucellier inversor.}
\label{Fig5}
\end{figure}

The most famous functional linkage is the Peaucellier inversor 
(see \cite[page 273]{HCV} and \cite[page 156]{CR}) depicted on Figure 
\ref{Fig5} (with $a^2 -r^2 =t^2$). 

The vertex $F=v_1$ is the only {\em fixed vertex} of the inversor, $Z:= (0)$.  
According to the 19-th century work on linkages, the Peaucellier inversor 
is supposed to be the functional for the inversion $J_t(z)= t^2/\bar{z}$ 
with the center at zero and radius $t$. 

Unfortunately this 
is not true for our definition of functional linkage because of the 
degenerate realizations 
 $\phi$, $\psi$  with $\phi(B)= \phi(D)$ and $\psi(A)= \psi(C)$. Note 
that there is a 3-torus of degenerate realizations $\psi$ 
with $\psi(A)= \psi(C)$,  
so even the dimension of $C(\L, Z)$ is not correct for a functional linkage 
with $n=m=1$. 

\begin{figure}[tbh]
\leavevmode
\centerline{\epsfxsize=3in\epsfbox{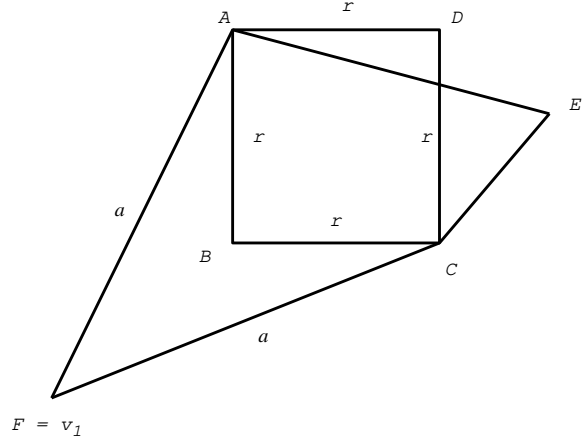}}
\caption{\sl The modified Peaucellier inversor $\J_t$: the square $ABCD$ is 
rigidified and $\ell[AE]- \ell[EC]= 2\eps >0$, $\ell[EC] > r$.}
\label{Fig6}
\end{figure}

Many of the degenerate realizations can be eliminated by rigidifying the 
square $ABCD$, but there remains an ${\mathbb S}^1\times {\mathbb S}^1$ of 
degenerate realizations 
with $\psi(A)= \psi(C)$ for which $\psi(B)$ and $\psi(D)$ are not in general 
related by inversion. We eliminate these by attaching a 
``hook''\footnote{Notice that by attaching this hook we had 
created an extra symmetry on the moduli space: the transformation 
which fixes images of all vertices except $\phi(E)$ and reflects 
$\phi(E)$ with respect to the line $(\phi(A)\phi(C))$.}  to $\{A, C\}$ 
as on the Figure \ref{Fig6}. 

\begin{lem}
\label{2.9}
The modified Peaucellier inversor $\J_t$ (with $B$ as input and $D$ as 
output) is a functional linkage for $J_t(z)= t^2/\bar{z}$. The domain 
of $\J_t$ is the annulus
$$
\{z\in \C: \rho\le |z|\le \rho^{-1}\}
$$
where $\rho= \sqrt{a^2 -\eps^2}-  \sqrt{r^2 -\eps^2}$. The only quasiwall 
of the inversor is the circle of inversion $\{z: |z|=t\}$. $Wall(\J_t)$ 
is the boundary of the annulus $Dom(\J_t)$. 
\end{lem}
\proof If the 
image of the square $[ABCD]$ under a realization $\phi$ is 
nondegenerate then $\phi(D)= J_t(\phi(B))$, see \cite{CR}, \cite{HCV}. Notice 
that $\phi(A)\ne \phi(C)$ for any realization $\phi$ (the hook!). 
If $\phi(B)= \phi(D)$ then we still have   
$\phi(D)= J_t(\phi(B))$ since for such a realization the point $\phi(B)$ 
bisects the segment $[\phi(A),\phi(C)]$. It is clear from the triangle 
inequalities that the above annulus equals $Dom(\J_t)$. 
Suppose that we are given $\phi(B)= z\in \C- \{0\}$.  
Then the location of $\phi(A), \phi(C)$ is uniquely determined up to 
reflection $\si_1$ in the line $\la$ through $0, z$ which interchanges 
these points. This determines the point  $\phi(D)\in \la$ as well. The 
reflection $\si_1$ determines an element of order 2 in $Sym(\J_t)$ which 
we denote by $\si_1$ as well. However the action of $\si_1$ on $C(\J_t, Z)$ 
is {\em free} since $\phi(A)\ne \phi(C)$ for all realizations $\phi$. 
Once the images of $A, B, C, D$ are determined there remains a single 
indeterminacy for $\phi$: we can reflect $\phi(E)$ via the reflection 
$\si_2$ in the line through $\phi(A), \phi(C)$. This reflection 
determines an element $\si_2\in Sym(\J_t)$, the fixed-point set of which 
consists of realizations for which $|\phi(A)- \phi(C)|$ is minimal, i.e. 
equal to $2\eps$. The projection of the set of these realizations to $Dom(\J_t)$ 
forms the boundary of the annulus $Dom(\J_t)$. This identifies the union of 
walls in $C(\J_t, Z)$. We also conclude that the  group of symmetries $Sym(\J_t)$ 
is $\Z_2\times \Z_2$ and is generated by the above involutions $\si_1, \si_2$. 

It remains to identify the quasiwalls. Suppose that 
$\xi\in T_{\phi}C(\J_t, Z)$ is an infinitesimal deformation which vanishes on 
$B$, recall that $\phi|v_1=0$ as well and $\phi$ does not belong to a wall. 
The above description of the domain of $\J_t$ implies that the triangles 
$$
\Del(0, \phi(B), \phi(A)), \quad \Del(0, \phi(B), \phi(C))
$$
are nondegenerate for each realization $\phi$ as above. Hence $\xi|A=0, \xi|C=0$. 
If the rhombus 
$$
[\phi(B), \phi(A), \phi(D), \phi(C)]
$$
is nondegenerate then Lemma \ref{tri} implies that $\xi|D=0$. Since 
$\phi$ is not on a wall, the triangle 
$$
\Del(\phi(E), \phi(C), \phi(A))
$$
is nondegenerate and $\xi|E=0$. This proves that the projection of the 
quasiwall to $Dom(\J_t)$ is the circle
$$
\{ \phi(D)= \phi(B)\} = \{ \phi(B): |\phi(B)|= t\}
$$ 
i.e. the circle of the inversion $J_t$. $\qed$

\begin{rem}
\label{inve}
Notice that unlike the cases of other linkages, the quasiwall $\{|z|=t\}$ 
does not move if we alter edge-lengths of the functional  linkage $\J_t$ for 
the given function $J_t$. On the other hand, by adjusting the parameters 
$\eps, a, r, \ell[AE]$ and keeping $t$ fixed we can get any point $z\in \C^*$ to the 
interior of $Dom(\J_t)$. 
\end{rem}

\begin{Not}
We shall use the notation $\J$ for $\J_1$. 
\end{Not}

\subsection{The multiplier}

Our construction of the multiplier is quite different from the one that 
was used by Kempe \cite{Kempe} and other people (see for instance \cite{B}, 
\cite{JS}). The idea is to use algebra instead of geometry: if one has addition, 
subtraction and inversion then one also gets the function $z\mapsto z^2$ via 
composition as follows. Consider the identity
$$
\frac{1}{\bar{z}-0.5} + [- \frac{1}{\bar{z}+0.5}]= \frac{1}{\bar{z}^2-0.25}
$$
Hence we can combine the following linkages: 

\begin{itemize}
\item Three translators for the functions $\tau_{\pm0.5}': z\mapsto z\mp 0.5$,  
$\tau_{0.25}: z\mapsto z+0.25$. 
\item One pantograph for the germ of the function $w\mapsto -w$ at 
the point $-2$. 
\item  Three inversors $\J$ for the function $J_1: z\mapsto 1/\bar{z}$. 
\item The modified adder $\L_A$ for the germ of addition at $(-2,-2)$. 
\end{itemize}

\noindent to get a functional linkage $\q$ for the function $z\mapsto z^2$. 

\begin{lem}
Under the following restrictions:
$$
\eps < r/2 , \quad a+r > \sqrt{3}/2 %8/\sqrt{3}
$$
on the parameters $a, r, \eps$ for the inversor $\J$,  
the origin belongs to $Dom^*(\q)$. 
\end{lem}
\proof First of all we need:
$$
\pm 0.5 \in Dom^*(\J)
$$
None of the points $\pm 0.5$ belong to the circle of inversion, hence 
if we use $a, r, \eps$ as above then 
$$
1> 0.5 > \rho 
$$
and $\pm 0.5\in Dom^*(\L)$. 

\medskip
The point $(-2, -2)$ belongs to $Dom^*(\L_A)$ provided that the parameter  
$t$ is chosen so that $t < 2 < 3t$. Finally, we apply the inversor 
$\J$ again to compute
$$
\frac{1}{\bar{z}^2-0.25}\mapsto z^2-0.25
$$
For this operation we need: $-4\in Dom^*(\J)$. Direct computation again 
shows that $4 < \rho^{-1}$ under the above restrictions on $a, r, \eps$,  
which implies that $-4\in Dom^*(\J)$. $\qed$ 

Thus, we have a functional linkage $\q$ for the computation of 
the function $z\mapsto z^2$ so that $0\in Dom^*(\L)$. Then we 
use the identity 
$$
zw= [(z+w)^2 +(- (z^2 + w^2)) ]/2
$$  
to construct a functional linkage for the complex multiplication. 
The linkages which are used for this computation are: the three 
copies of modified adder, two modified pantographs (for the functions 
$x\mapsto -x$ and $y\mapsto y/2$) and three copies of the linkage 
$\q$ for squaring. On each step of the composition of linkages all 
we need functional linkages with the origin in $Dom^*$, which is true 
for all the above linkages.   

\subsection{The straight-line motion linkage}
\label{stra}

In this section we modify the usual 
{\em Peaucellier straight-line motion linkage} (see \cite{CR}, \cite{HCV}) 
to obtain a real functional linkage ${\cal S}$ for the inclusion 
$\R \to \C = \R^2$.  
  
\begin{figure}[tbh]
\leavevmode
\centerline{\epsfxsize=3in\epsfbox{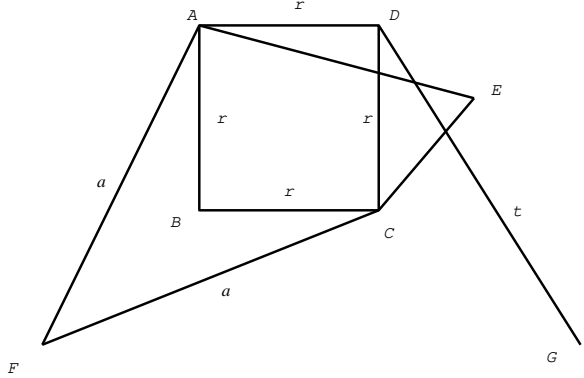}}
\caption{\sl The Peaucellier straight-line motion linkage ${\cal S}$. 
The vertices $F, G$ are fixed. The vertex $B$ is the input and output.} 
\label{Fig7a}
\end{figure}

\medskip
Start with the rigidified inversor $\J_t$ and add the edge $[GD]$ of the 
length $t$. The vertices $F, G$ are the {\em fixed vertices} of the new 
linkage ${\cal S}$.  Their images are: 
$\phi(F)= - \phi(G)= \pm \sqrt{-1} t/2$. Take the 
vertex $B$ as both the input vertex and the output vertex. See Figure 
\ref{Fig7a}.

\begin{rem}
This choice is somewhat strange from the classical point of view since 
the linkage ${\cal S}$ was invented to transform periodic linear 
motion of the vertex $B$ to the circular motion of the vertex $D$ 
(from this point of view $B$ is the input and $D$ is the output). 
However we do not 
use the linkage ${\cal S}$ to transform linear to circular 
 motion but to restrict motion of the input-vertex $B$ to the real axis. 
\end{rem}

The point $\phi(D)$ is now restricted 
to the circle with the center at $\phi(G)$ and radius 
$t$. The input $\phi(B)$ 
is obtained from $\phi(D)$ by inversion with the center at $\phi(F)$ and 
radius $t$. Whence the input $\phi(B)$ moves along a segment in the real axis. 
Notice that ${\cal S}$ has a symmetry that interchanges $\phi(F)$ and $\phi(G)$, 
and maps $\phi(D)$ to $\overline{\phi(D)}$; this symmetry is induced by 
 the complex conjugation $\C \to \C$. 

We use the following restrictions on the side-lengths of the linkage:
$$
0< 2\eps = \ell[AE]- \ell[CE]\ ,$$ 
$$
\ell[CE] > 2 r, \quad 
a> r > \eps, \quad 17r > 15 a
$$
Under these conditions the linkage ${\cal S}$ is a real functional linkage 
for the inclusion map $id:\R \to \C=\R^2$ and 
the input map $p: C({\cal S},Z) \to \R\subset \R^2$ has the following property: 

\smallskip
$Dom^*({\cal S},Z)$ contains the open interval 
$(-\frac{\sqrt{3}}{2}t, \frac{\sqrt{3}}{2}t)$. 

\begin{figure}[tbh]
\leavevmode
\centerline{\epsfxsize=3in\epsfbox{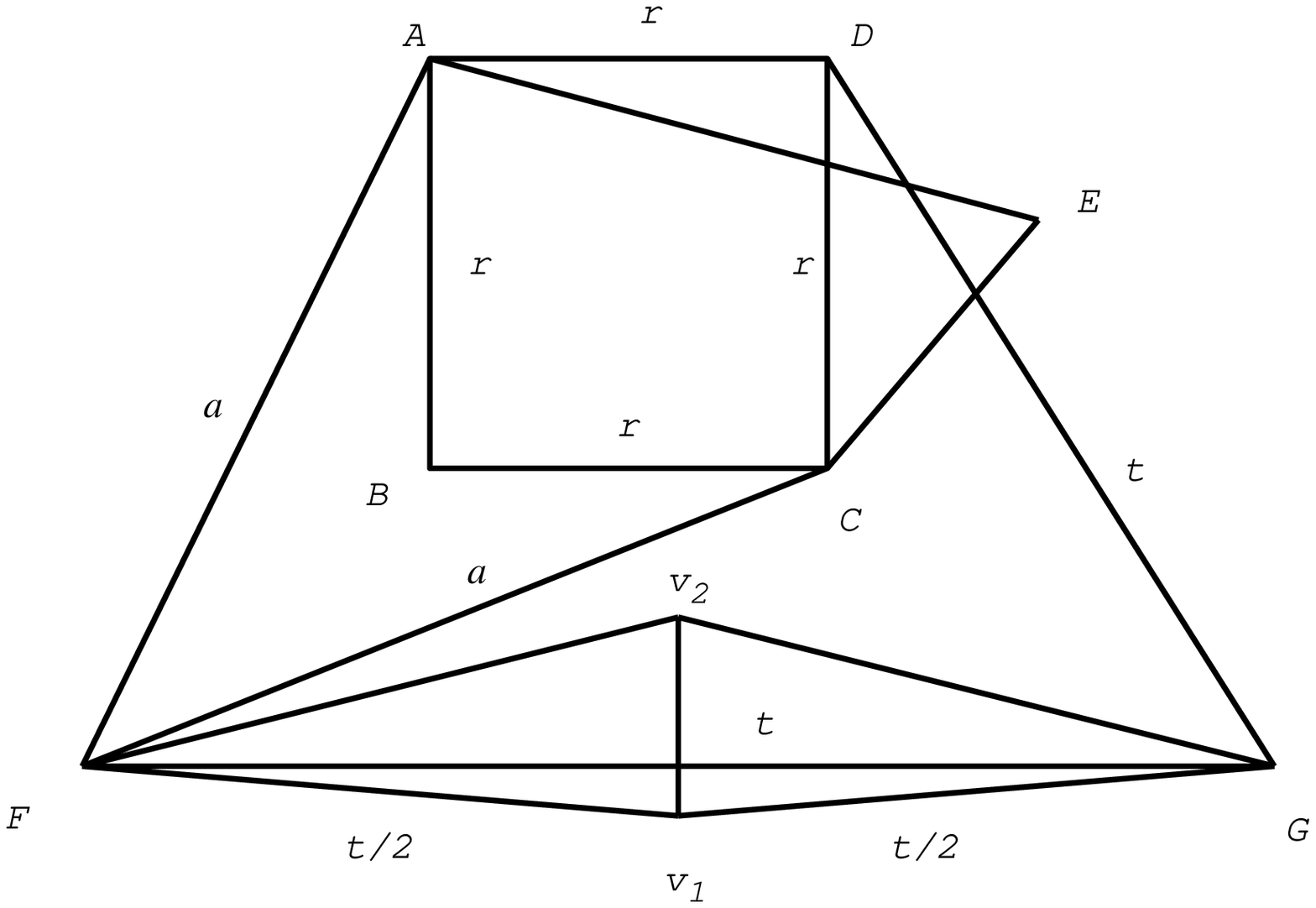}}
\caption{\sl The based Peaucellier straight-line motion linkage $\t{\cal S}$:  
$\ell[v_1 v_2]^2 + (t/2)^2 = \ell[Fv_2]^2= \ell[Gv_2]^2$. 
The image of $B$ under all realizations lies on a 
segment of the real axis which contains the open interval 
$(-\frac{\sqrt{3}}{2}t, \frac{\sqrt{3}}{2}t)\subset Dom^*(\t{\cal S})$.}
\label{Fig7}
\end{figure}

\medskip
Notice that $\phi(F), \phi(G)\notin \R$. We will need a modification 
$\t{\cal S}$ of ${\cal S}$ where images of all fixed vertices are real 
numbers. The {\em based} linkage $\t{\cal S}$ is produced from ${\cal S}$ 
via the construction in Section \ref{fixing}, see Figure \ref{Fig7}. 

We will also need another modification ${\cal S}^m$ of the 
linkage $\t{\cal S}$. Namely, take $m$ isomorphic copies $\t{\cal S}_j$ 
of the linkage $\t{\cal S}$. Then take their fiber sum by identifying 
the vertices with the labels $v_1, v_2, F, G$ for each pair $\t{\cal S}_j$, 
$\t{\cal S}_i$. We leave it to the reader to verify that
 ${\cal S}^m$ is a real functional linkage for the inclusion map 
$id: \R^m \to \C^m$ and
$$
Dom^*({\cal S}^m)\supset (-\frac{\sqrt{3}}{2}t, \frac{\sqrt{3}}{2}t)^m
$$

\noindent The linkage ${\cal S}^m$ is used for constructing {\em real} 
functional linkages from the {\em complex} ones.  

\section{Expansion of domains of functional linkages}

\begin{conv}
In this section we will suppress choices of fixed vertices and their images for 
functional linkages.  
\end{conv}

We apply the results of Section \ref{pantograph} 
to expand domains of functional linkages:

\begin{lem}
Suppose that $g(x)$ is a homogeneous polynomial of degree $d$, 
$\L$  is a functional linkage 
which defines the germ $(g, 0)$. Then for any $r>0$ we can modify $\L$ 
so that the new linkage $\t\L$ is functional for the function $g$ and 
$Dom^*(\t\L)$ contains the disk $B_r(0)$. 
\end{lem}
\proof We consider the case $d>0$, the case of constant functions 
($d=0$) is left to the reader. 
By the assumption $Dom^*(\L)$ contains a disk $B_{\eps}(0)$ 
centered at the origin, we can assume $\eps< r$. 
Choose positive $\la< \eps/r <1$. Let $\mu:= \la^{-d} > 1$. 
We use the formula 
$$
g(y)= \la^{-d} g(\la y)= \mu g(\la y)
$$
to construct a functional linkage $\t\L$ for the function $g$ as a 
composition of the following linkages:

\begin{itemize}
\item $\p_{\la}'$ (the modified pantograph for the multiplication by $\la$), 
\item the linkage $\L$,  
\item $\p_{\mu}'$ (the modified pantograph for the multiplication by $\mu$). 
\end{itemize}

\noindent The linkages $\p_{\la}'$, $\p_{\mu}'$ are chosen so that 
$$
B_r(0)\subset Dom^*( \p_{\la}'), \quad B_R(0)\subset Dom^*(\p_{\mu}' )
$$ 
where $R:= \max \{ g(y): y\in B_{\eps}(0)\}$. 
Let's check that this choice guarantees that $B_r(0)$ is contained in 
$Dom^*(\t\L)$. By Theorem \ref{funthm1} we have:
$$
Dom^*(\t\L) \supset [Dom^*( \p_{\la}')] \cap [\la^{-1}(Dom^*(\L))] \cap 
[\la^{-1} g^{-1}(Dom^*(\p_{\mu}' ))] 
$$
Then, since $\la^{-1}\eps > r$, it is enough to verify that
$$
B_r(0)\subset \la^{-1} g^{-1}(Dom^*(\p_{\mu}' ))
$$ 
However $B_R(0)\subset Dom^*(\p_{\mu}' )$, hence 
(by the choice of $R$)
$$
B_{\eps/\la}(0)\subset \la^{-1} g^{-1}(Dom^*(\p_{\mu}' ))
$$
which together with the inequality $r\le \la^{-1}\eps$ implies the 
assertion.  \qed  

As a corollary we get the following Theorem:

\begin{thm}
\label{exp}
(Theorem on expansion of domain.) 
Suppose that $f: \k^m \to \k^n$ be a polynomial morphism, 
$\L$  is a functional linkage which defines the germ $(f, 0)$. Then 
for any $r>0$ we can modify $\L$ so that the new linkage $\t\L$ is 
functional for the morphism $f$ and $Dom^*(\t\L)$ contains the disk 
$B_r(0)$. 
\end{thm}
\proof We will consider the case when $n=1$, the general case follows from 
Theorem \ref{funthm2}. Write $f(x)$ as 
$$
f(x)= \sum_{j\le d} f_j(x)$$ 
where each $f_j$ is a homogeneous polynomial of degree $j$. 
Let $g(y):= y_1 + ... + y_d$. Hence we can represent $f$ as a 
composition of homogeneous polynomials $f_j, j\le d,$ and $g$. Now 
the assertion follows from the previous lemma and Theorems 
\ref{funthm1}, \ref{funthm2}.  \qed  

\section{Realization of complex polynomial maps by functional linkages}
\label{complex}

In this section we prove Theorem A (the complex case).  
We first consider the case $f: \C^m \to \C$, i.e. $n=1$. Let
$$
f(x)= a_0+ \sum_{j} a_j g_j(x)
$$
where $g_j=x_1^{\al_1}... x_m^{\al_m}$ are monomials of 
positive degrees and $a_j\in \C$ are constants ($j=0, 1,...,N$). 
Let $y= (y_0,...,y_N)$. Consider the function 
$$
\hat{f}(x, y)= y_0+ \sum_{j} y_j g_j(x)
$$
This function is obtained via composition of the multiplication and 
addition operations. Hence we use the elementary linkages for the addition and 
multiplication we get a complex functional linkage  $\hat\L$ for the germ 
$(\hat{f}, 0)$. Then we use Theorem \ref{exp} (on expansion of domain): 
for each given $\rho >0$ 
we can modify $\hat\L$ to $\t\L$ so that $\t\L$ is functional for the pair 
$(\hat{f}, B_{\rho}(0))$, $0\in \C^{m+N}$. We use $\rho$ so large that 
$B_{\rho}(0)$ contains the disk
$$
\{(x,y): x\in B_r(\O), y_j=a_j, j=0,..., N\}
$$
We represent $f$ as a composition of the function $\hat{f}$ and 
the constant function
$$
a: (y_0,..., y_N)\mapsto (a_0,..., a_N)
$$
The constant function is defined by a functional linkage as follows:

Let $\a$ be the graph which consists of the set of vertices 
$[In(\a)= (P_1,..., P_m)]\cup [Out(\a)= (Q_1,..., Q_N)]$, no edges, 
$W= Out(\a)$ and $Z= (a_0,..., a_N)$.   

Clearly $Dom^*(\a)= \C^m$. Thus the 1-st and 2-nd functionality 
theorems \ref{funthm1}, \ref{funthm2} imply that composition of the linkages 
$\t\L$ and $\a$ gives us a functional linkage for the pair $(f, B_r(\O))$. 
This proves Theorem A in the complex case when $n=1$. 

To get functional linkages for polynomial vector-functions we 
use repeatedly the 2-nd functionality theorem \ref{funthm2}: 

If we have a functional linkage $\L_1$ for the germ 
$(f_1(x_1,..., x_n), A_1)$ and a functional linkage $\L_2$  for 
$(f_2(x_1,..., x_n), A_2)$ we glue inputs of $\L_1$ and $\L_2$ to construct 
a functional linkage $\L$ for the germ  $((f_1,f_2), (A_1, A_2))$. 

Thus we proved

\begin{thm} 
\label{mainc}
Let $f: \C^n \to \C^m$ be a polynomial map, $\O\in \C^m$ and $r>0$. 
Then there is a marked functional linkage $\L= (L, \ell, W)$ together 
with a vector $Z\in \C^s$ so that: 

The ball $B_r(\O)$ is contained in $Dom^*(\L, Z)$, 
$q\circ p^{-1}: Dom(\L, Z)\to \C^m$ equals the restriction  
of the vector-function $f$, i.e. $(\L, Z)$ defines $(f, B_r(\O))$. 
\end{thm} 

There is a special case when $f$ has real coefficients and $Z\in \R^s$. 
Recall that we use only real numbers $b$ for the 
translators in $\L$. We apply the construction described in the 
Section \ref{fixing} and modify $\L$ to a based linkage $\t\L$. According to 
Lemma \ref{fi} we get an algebraic isomorphism 
$$
\tau: \M(\t\L)\to C(\L, Z)
$$
Hence the based linkage $\t\L$ also defines the pair $(f, B_r(\O))$. 
This proves Theorem A of the Introduction (in the complex case).

\section{Transition from complex to real functional linkages}
\label{real}

Let $f: \R^m \to \R^n$ be a polynomial function, 
$\O\in \R^m$ be a point and $r>0$. 
Our goal is to produce a real functional linkage for the 
pair $(f, B_r(\O))$. 
 
We extend $f$ by complexification to a morphism $f^c: \C^m \to \C^n$ and 
construct a complex functional linkage $\L'$ for $(f^c, B_r^c(\O))$, where  
$B_r^c(\O)$ is the ball of radius $r$ in $\C^m$ centered at $\O$. 

Take the real functional linkage ${\cal S}^m$ for the identity 
map $id: \R^m \to \R^m$, see \S \ref{stra}. 
We choose the parameter $t$ in ${\cal S}^m$ so large that 
$$
B_r(\O)\subset (-\frac{\sqrt{3}}{2}t, \frac{\sqrt{3}}{2}t)^m 
\subset Dom^*({\cal S}^m)
$$ 
Next we alter $\L'$ via fiber product with the linkage ${\cal S}^m$. 
Namely, take the bijection 
$\be: In(\L')\to In({\cal S}^m)= Out({\cal S}^m)$ which 
maps each input vertex $P_j'$ of $\L'$ to the input vertex $P_j''$ of 
${\cal S}^m$. Let $\L:=  \L' *_{\be} {\cal S}^m$. Then the 1-st 
functionality theorem \ref{funthm1}  implies that $\L$ is a real-functional 
linkage for the polynomial $f$ and $Dom^*(\L)$ contains the disk 
$B_r(\O)\subset \R^m$. 

Notice that for each fixed vertex $v\in W$ of ${\cal S}^m$ we have: 
$z=\phi(v)\in \R$, $\phi\in C({\cal S}^m, Z)$. The same is 
true for $\L'$ since the polynomial $f$ has only real coefficients 
and we use only real numbers $b$ for the translators. Then (as in as in 
the Section \ref{complex}) Lemma \ref{fi} implies that we can modify 
$\L$ to a based linkage $\t\L$ so that
$$
\tau: \M(\t\L)\to C(\L,Z)
$$
is an algebraic isomorphism. This concludes the proof of 
Theorem A in the real case.   \qed  

\section{Realization of algebraic sets and 
smooth manifolds as moduli spaces of planar linkages}
\label{Kemp}

In this section we derive Theorem B from Theorem A.  
Let $M$ be a compact affine algebraic subset of $\R^m$, choose 
a polynomial $f: \R^m \to \R$ such that $M= f^{-1}(0)$. We may 
assume $M\subset B_r(0)$. 
By Theorem A (the real case) we have a {\em based} functional 
linkage $\L$ for the pair $(f, B_r(0))$. The output mapping $p$ of $\L$ is 
an analytically trivial polynomial covering over $B_r(0)$. 
Now glue the output vertices of $\L$ to the basic vertex $v_1$ to 
obtain a linkage $\L_0$. Let $p_0$ be the output mapping of $\L_0$. 
The images of the input vertices of $\L_0$ 
under $\phi\in \M(\L_0)$ are now constrained to $f^{-1}(0)$. 
The mapping $p_0: \M(\L_0) \to \R^m$ is an analytically trivial 
polynomial covering  over $M$ (by the 3-rd functionality theorem \ref{funthm3}). 
 Theorem B follows. $\qed$ 

\medskip
Proof of Theorem B' is similar and is left to the reader.

\bigskip
To prove Corollary C we use Theorem \ref{nash} and Theorem B to get a 
linkage $\L_0$ and an analytically trivial covering $p: \M(\L_0) \to M$ 
where we identify $M$ with a real algebraic subset of $\R^m$. $\qed$ 

\section{How to draw algebraic curves}
\label{draw}

In this section we prove Theorem \ref{kem} according which one can 
``draw'' arbitrary algebraic curves in $\R^2$ using planar 
mechanical linkages. For instance, if $\Ga$ is a compact connected 
algebraic curve in $\R^2$ then there is a closed complex functional 
linkage $\L_0$ with a single input vertex $P$ so that as realizations 
$\psi$ of $\L_0$ vary along an arbitrary connected 
component $C$ of $\M(\L_0)$, the vertex $P$ traces the curve $\Ga$ and 
the projection $C\to \Ga$ is an analytic isomorphism.

We first need a functional linkage for the complex conjugation. 
There are several ways to do it. 

\medskip
{\bf The 1-st construction.} Using Theorem B' construct a closed functional 
linkage $\L_0$ for the germ of the complex algebraic set $zw=1$ at  the point 
$(2, 1/2)\in \C^2$. Let 
$P_1, P_2$ be the input vertices of $\L_0$. Then let $\L^0$ be the linkage 
$\L_0$ where we declare $P_1$ the sole input vertex and $Q_1:=P_2$ the output 
vertex. Then $\L^0$ is a complex functional linkage for the germ of the 
function $z\mapsto z^{-1}$ at the point $2$. Recall that we have the complex 
functional linkage $\J$ for the germ of the inversion
$$
w\mapsto 1/\bar{w}
$$
at the point $1/2$. Hence we compose 
$z\mapsto z^{-1} \mapsto \bar{z}$ and compose the linkages  
$\J$, $\L^0$ to get a $\C$-functional linkage $\bar{\L}$ 
for the germ of 
the map $z\mapsto \bar z$ at  the point $2$. Finally, we use the formula
$$
\bar{z}= \overline{z+2} - 2
$$ 
and composition of $\bar{\L}$ with two translators to get a $\C$-functional 
linkage $\L'$ for the germ $(z\mapsto \bar{z}, 0)$. 

\begin{figure}[tbh]
\leavevmode
\centerline{\epsfxsize=3in\epsfbox{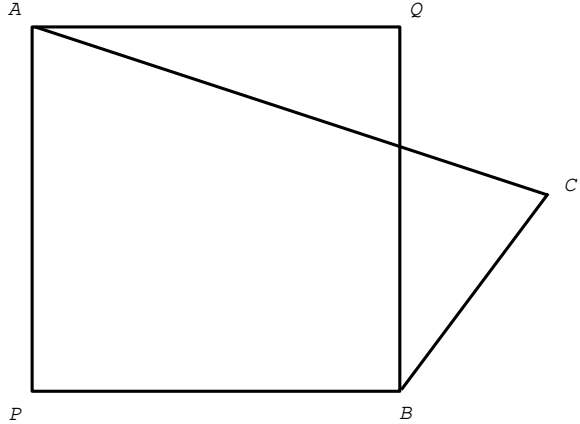}}
\caption{\sl Construction of a functional linkage for the complex conjugation.}
\label{F21}
\end{figure}

\bigskip
{\bf The 2-nd construction.} We start with the linkage $\b$ described on Figure 
\ref{F21}: an abstract (rigidified) square with a ``hook'' attached. We let 
$$
\ell[AP]= \sqrt{2}, \quad |\ell[AC]-\ell[BC]| > 2$$ 
Let ${\cal S}^2$ be the linkage from the section \ref{stra} where $t=1$, 
let $P_1, P_2$ be the input vertices of ${\cal S}^2$ (they are the output 
vertices as well). Take $\be: A\mapsto P_1, B\mapsto P_2$ and 
$\L:= \b *_{\be} {\cal S}^2$. We declare $P\in \b$ the input and $Q\in \b$ 
the output of the linkage $\L$. We leave it to the reader to verify 
that the linkage $\L$ is a $\C$-functional linkage for the germ of 
$z\mapsto \bar z$ at the point $\sqrt{-1}$.

\medskip

\begin{rem}
Under all realizations $\phi$ of $\L$ we have: $\phi(A), \phi(B)\in \R$, 
$\phi(A)\ne \phi(B)$ and $\phi(Q)= \overline{\phi(Q)}$.  
\end{rem}

Finally, we use the formula
$$
\bar{z}= \overline{z+i} - i
$$ 
and composition of $\L$ with two translators to get a $\C$-functional linkage 
$\L''$ for the germ $(z\mapsto \bar{z}, 0)$. 

\begin{thm}
\label{kem}
Let $f= f(z, \bar{z}), f:\C \to \R$ be a polynomial function of the variables 
$z, \bar{z}$ and $\Ga:= f^{-1}(0)\subset \C$ be a real-algebraic curve. Pick an 
open (in the classical topology) bounded subset $U\subset \Ga$. Then there is a 
closed $\C$-functional linkage $\L_0$ so that the input map  
$p_0: C(\L_0, Z)\to \C$ is an analytically trivial polynomial covering over 
$U$. 
\end{thm}

\proof Let $U\subset B_r(\O)$. Our argument is exactly the same as in the 
proof of Theorem B. Namely, as in Theorem A we first construct a functional 
linkage $\L$ 
for the pair $(f, B_r(\O))$ (now we use the composition of addition, multiplication 
and the complex conjugation). Then we attach the output vertex $Q$  of $\L$ to the 
distinguished vertex $v_1$ (such that $\phi(v_1)=0$ for all relative realizations). 
Let $\L_0$ be the resulting closed functional linkage and $p_0$ be its input map. 
Then as in the proof of Theorem B we have: $p_0$ is an analytically trivial 
polynomial covering over $U$.  \qed

\section{Universality theorem for arrangements in $\P^2$}
\label{arrangements}

In the section we review notions of configuration spaces for arrangements 
and universality theorems proven in  \cite{KM6}, which extend earlier results 
of Mnev \cite{Mnev}. 

%\section{The Moduli Space of a Planar Arrangement}

Let $\a$ be an {\em abstract arrangement}, i.e. a
 bipartite graph with the parts $\p$ and $\L$. We say that a ``point'' 
$P\in \p$ is incident to a ``line'' $L\in \L$ if $P$ and $L$ are 
connected by an edge. A projective realization $\phi$ of $\a$ is a map
$$
\phi: \p \cup \L \to \P^2 \cup (\P^2)^{\vee},\ \ \phi(\p)\subset \P^2,\ \  
\phi(\L)\subset (\P^2)^{\vee}
$$
such that if $P$ and $L$ are incident then $\phi(P)\in \phi(L)$. 
This condition defines a projective scheme $R(\a)$ 
over $\Z$. We let $R(\a, \C\P^2)$ and $R(\a, \R\P^2)$ denote the sets of 
complex and real points of $R(\a)$. 

Here and below we use the symbol $\v$ for polarity between points and lines 
in the projective plane, this polarity is determined by the standard bilinear  
form on $\R^3$ given by $\|(x,y,z)\|^2= x^2 +y^2 +z^2$.  

%We will also use the term {\em projective arrangements} for projective 
%realizations.  

We now want to pass to the quotient of $R(\a)$ by $PGL_3$. We do this by 
restricting to realizations in a ``general position'' and then taking a 
cross-section. To make it precise we first define based arrangements.

\begin{figure}[tbh]
\leavevmode
\centerline{\epsfxsize=6in\epsfbox{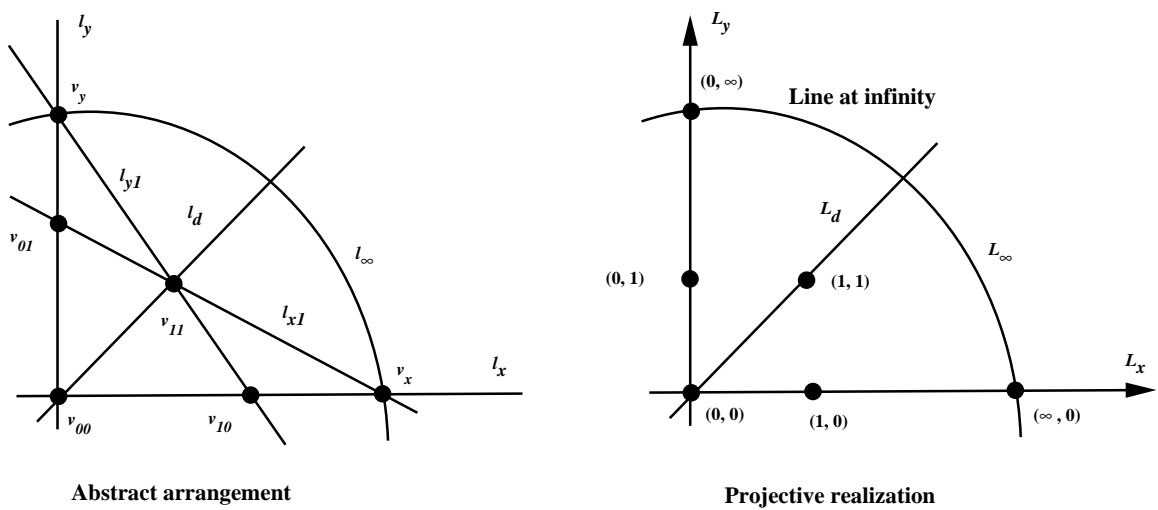}}
\caption{\sl The standard triangle $T$ and its standard realization.}
\label{Fig23}
\end{figure}

\begin{defn}
\label{5.1}
The {\bf standard triangle} is the abstract arrangement $T$ consisting of 
6 point-vertices and 6 line-vertices that corresponds to a triangle with its 
medians, see Figure \ref{Fig23}.  
\end{defn}

\begin{defn}
\label{5.2}
The {\bf standard realization} $\phi_T$ of the standard triangle $T$ is 
determined by:
$$
\phi_T(v_{00})= (0,0), \phi_T(v_{x})= (\8,0), \phi_T(v_{y})= (0,\8), 
\phi_T(v_{11})= (1,1)
$$
\end{defn}

Here $(0,0), (\8,0), (0,\8), (1,1)$ are points in the affine plane 
$\A^2\subset \P^2$  which have the homogeneous coordinates: $(0:0:1), 
(1:0:0), (0: 1: 0), (1:1:1)$ respectively. 

We say that an abstract arrangement $\a$ is {\em based} if it comes 
equipped with an embedding $i: T\to \a$. Let $(\a, i)$ be a based 
arrangement. We say that a projective realization $\phi$ of $\a$ is 
{\em based} if $\phi\circ i= \phi_T$. Let $BR(\a, \P^2(\k))$ be the 
subset of $R(\a, \P^2(\k))$ consisting of based realizations, $\k=\R, \C$. 

\begin{lem}
(See \cite[Theorem 8.20]{KM6}.) $BR(\a, \P^2(\C))$ is the set of complex points 
of a projective scheme over $\Z$ which is a scheme-theoretic quotient of $R(\a)$ 
by the action of $PGL_3$. 
\end{lem}

\begin{defn}
\label{5.3}
 A {\bf functional arrangement} is a based arrangement $({\a},i)$ with 
two subsets of marked point-vertices $\mu= (P_1,..., P_m)$ (the {\bf 
input-vertices}) and point-vertices 
$\nu= (Q_1,..., Q_n)$ (the {\bf output-vertices}) 
such that all the marked vertices are incident 
to the line-vertex  $l_x \in i(T)$ (which corresponds to the $x$-axis) 
and such that the following two axioms are satisfied:

Let $BR_0(\a_{\mu})\subset BR(\a)$ denote the 
open subset which consists of realizations $\phi$ such that 
$\phi(P_j)\in \A^2$ for all $j$, we define $BR_0(\a_{\nu})$ 
similarly. Then we require: 

(1) $BR_0(\a_{\mu}) \subset BR_0({\a}_{\nu})$.

(2) The projection $p: BR_0({\a}_{\mu}) \to \A^m$ given by $p(\phi)= 
(\phi(P_1),..., \phi(P_m))$ is an isomorphism of schemes over $\Z$. 
\end{defn}

Each functional arrangement determines a morphism $f: \A^m \to \A^n$ 
(which is defined over $\Z$) by the formula:
$$
f(x)= q\circ p^{-1}(x)
$$ 
where $q(\phi)= (\phi(Q_1),..., \phi(Q_n))$.

\begin{thm}
(See \cite[Lemma 9.7]{KM6}.) Let $f: \A^m \to \A^n$ be any polynomial 
mapping with integer coefficients. Then there is a functional arrangement 
$\a$ which determines $f$.
\end{thm}

\bigskip
Let $S\subset \A^m$ be a closed subscheme defined over $\Z$, $S= f^{-1}(0)$ 
for some morphism $f: \A^m\to \A^n$. Let $\a$ be a functional 
arrangement which determines $f$ as in the above theorem. 
By gluing the output vertices of $\a$ to $v_{00}$ we obtain an arrangement 
$\a^0$ containing distinguished vertices $P_1,..., P_m$. Again define 
$BR_0(\a^0)$ by requiring $\phi(P_i)$ to be finite. We get an induced 
morphism (easily seen to be an embedding) $p: BR_0(\a^0)\to \A^m$. 
We then have

\begin{thm} 
\label{moremnev}
(Theorem 1.3 of \cite{KM6}) Let $S$ be a closed subscheme of $\A^m$ 
(again over $\Z$). Then there exists a based marked arrangement $\a$ such 
that the input mapping $p: BR_0(\a)\to \A^m$ induces an isomorphism 
of schemes $BR_0(\a)\to S$.
\end{thm}

We will use the following version of the above theorem:  

\begin{thm}
\label{F}
Let $X$ be a compact real algebraic set 
defined over $\Z$. Then there exists a based arrangement $\a$ 
such that $X$ is entire birationally isomorphic to a Zariski open and 
closed subset $C$ in $BR(\a, \P^2(\R))$. 
\end{thm}
\proof Using Theorem \ref{closed} we may assume that $X$ is 
projectively closed. We choose the projective scheme 
${\mathfrak X}\subset \P^m$  whose set of real points is $X$ 
so that the corresponding affine scheme 
${\mathfrak X}_a\subset  \A^{m+1}$ corresponds to a real reduced ideal. Thus $X$ is Zariski dense in ${\mathfrak X}(\C)$. 

Define the affine scheme 
${\mathfrak X}_A={\mathfrak X}\cap \A^m$.  
Now apply Theorem \ref{moremnev} to construct a 
based marked arrangement $\a$ so that  $BR_0(\a)$ is isomorphic to 
${\mathfrak X}_A$ (as a scheme), hence the sets of real points of 
these schemes are isomorphic as well. Thus $X$ is (polynomially)  
isomorphic to $BR_0(\a, P^2(\R))$ which is Zariski open. It remains 
to show that $BR_0(\a, P^2(\R))$ is also Zariski closed. 

Recall that $BR(\a)$ embeds canonically in a product 
$(\P^2)^N \times (\P^1)^m$ where the last $m$ factors correspond to 
the input vertices. The morphism $p: BR_0(\a)\to \A^m$ is the restriction 
of the projection on the last $m$ factors:  
$$
(\P^2)^N \times (\P^1)^m \to (\P^1)^m 
$$
The subset $BR_0(\a, \P^2(\C))$ is constructible, hence its closure with 
respect  to the classical topology is the same as its closure 
$\overline{BR_0}(\a, \P^2(\C))$ with 
respect to the Zariski topology in $(\P^2(\C))^N \times (\P^1(\C))^m$. 

\noindent Suppose that there is a real point 
$z\in \overline{BR_0}(\a, \P^2)$  
that does not belong to $BR_0(\a, \P^2(\R))$. Then $z$ is the limit of  
a sequence  $z_j\in BR_0(\a, \P^2(\C))$. However $p(z_j)\in \C^m$ are 
obtained by ``forgetting'' all but the last $m$ coordinates of $z_j$, hence 
$p(z_j)$ will converge to a real point $x$ of ${\mathfrak X}$. It is clear that 
$x\notin \R^m$ and hence does not belong to $X$. This contradicts the fact 
that $X$ is projectively closed.  \qed  

\begin{rem}
In general  $BR(\a, \P^2(\C))$  is different from 
$\overline{BR_0}(\a, \P^2(\C))$. As an example consider the linkage 
$\a$ corresponding (via the construction in \cite{KM6}) to the system of 
equations:
$$
x+y =0, \quad x+y =1, \quad x =y 
$$
in $\A^2$. The set solutions of this system of equations is empty 
(even in the projective compactification of $\A^2$). 
Thus $BR_0(\a, \P^2(\C))= \0$, on the other  hand: 
$BR(\a, \P^2(\C))$ is a single point. 
\end{rem}

\bigskip
Now we construct  metric graphs corresponding to based abstract arrangements. 
Suppose that $\a$ is a based arrangement. We start by identifying the  
point-vertex  $v_{00}$ with the line-vertex  
$l_{\8}$, the point-vertex $v_x$ with the line-vertex $l_y$ and the 
point-vertex $v_y$ with the line-vertex  $l_x$ in the standard triangle 
$T$. We also introduce  the new edges 
$$
[v_{10} v_{00}], \quad [v_{01} v_{00}], \quad 
[v_{10} v_{x}], \quad [v_{01} v_{y}]
$$
(Here $v_{10}, v_{00}, v_{11}, v_{01}, ...$ 
are the {\em point}-vertices in the 
standard triangle $T$.) We will use  the notation $L$ for the 
resulting graph. We construct a length-function $\ell$ on the set of 
edges $e\subset L$ as follows: 

1) We assign the length $\pi/4$ to the new edges. 

2) We assign the length $\pi/2$ to the rest of the edges.

\section{Relation between two universality theorems}
\label{88}

The goal of this section is to establish a relation between two 
universality theorems for realizability of real algebraic sets 
(Theorems B and \ref{F}). 
Consider an abstract based arrangement $\a$. 
We choose $v_{00}, v_{x}, v_y, v_{01}, v_{10}$ 
as distinguished vertices of the 
corresponding metric graph $L$. Let $\L$ denote the metric graph $L$ 
with the distinguished set of vertices as above. 
Let $X$ be either ${\mathbb S}^2$ or $\R\P^2$ 
with the standard metric $d$ (so that the standard projection 
${\mathbb S}^2\to \R\P^2$ is a local isometry). 
Define the {\em configuration space}  
$C(\L, X)$ of realizations  of $\L$ in $X$ to 
be the collection of mappings $\psi$ from the 
vertex-set $\V(\L)$ of $\L$ to $X$ such that
$$
d(\psi(v), \psi(w))^2= (\ell[vw])^2
$$
for all vertices $v, w$ of $\L$ connected by an edge. 

\begin{rem}
Notice that if $a, b\in \R\P^2$ are within the distance $\pi/2$ then 
there are two minimal geodesics connecting $a$ to $b$. This is the 
reason to define $C(\L, X)$ as the set of maps from $\V(\L)$ rather 
than from $\L$ itself. 
\end{rem}

One can easily see that $C(\L, X)$ has natural structure of a real 
algebraic set. The subsets 
$$
\M(\L, \R\P^2):= \{\psi \in C(\L, \R\P^2): \psi(v_{00})=(0,0), 
\psi(v_{x})=(\8,0), $$
$$
\psi(v_{10})=  (1,0), \psi(v_{01})= (0,1)\}
$$
$$
\M(\L, {\mathbb S}^2):= \{\psi \in C(\L, {\mathbb S}^2): \psi(v_{00})=(0,0,1), 
\psi(v_{y})=(0, 1, 0), $$
$$
\psi(v_{x})=(1, 0, 0), 
\psi(v_{10})=  (1, 0, 1), \psi(v_{01})= (0, 1, 1)\}
$$
form cross-sections to the actions of the groups of isometries 
$PO(3,\R), O(3,\R)$ of $X$ on $C(\L, X)$. We call $\M(\L, X)$, 
 the {\em moduli spaces} of realizations of $\L$ in $X$ (where $X= 
{\mathbb S}^2, \R\P^2$). 

\begin{rem}
Now it is convenient to use the full group of isometries of ${\mathbb S}^2$ 
instead of the group of orientation-preserving isometries that we 
used for planar linkages. 
\end{rem}

\begin{lem}
The moduli space $\M(\L, \R\P^2)$ is (polynomially)  
isomorphic to the real algebraic set $BR(\a, \R\P^2)$.
\end{lem}
\proof The key to the proof is the fact that a point $P\in \R\P^2$ is 
incident to a line $L\in (\R\P^2)^{\v}$ iff 
$$
d(P, L^{\v})= \pi/2
$$
Thus we construct a morphism 
$$
\mu: BR(\a_0, \R\P^2) \to \M(\L, \R\P^2), \quad \mu: \phi\mapsto \psi
$$
so that for each point-vertex $P\in \a$ we have $\psi(P)= \phi(P)$ 
and for each line-vertex $L\in \a$  we have $\psi(L)= \phi(L)^{\v}$. 
This morphism has algebraic inverse given by the same formula 
(since $(L^{\v})^{\v}= L$).  \qed  

\medskip
Let $\M_0(\L, \R\P^2)$ be the image of $BR_0(\a_0, \R\P^2)$ 
under the isomorphism $\mu$. 
Consider the standard 2-fold covering ${\mathbb S}^2 \to \R\P^2$. 
It induces a (locally trivial) analytical covering 
$$
\al: \M(\L, {\mathbb S}^2)\to \M(\L, \R\P^2)
$$
The group of automorphisms of $\al$ is $(\Z_2)^r$, where $r$ is 
the number of (point) vertices in $[L - \p(T)]\cup \{v_{11}\}$. 
The generators of this group are indexed by the vertices $v\in 
[L - \p(T)]\cup \{v_{11}\}$: 
$$
g_v: \psi(v)\mapsto - \psi(v),  g_v: \psi(w) \mapsto \psi(w), w\ne v
$$

\begin{prop}
For each arrangement $\a$ as in Theorem \ref{F}, 
the covering $\al$ is analytically trivial over $\M_0(\L, \R\P^2)$. 
\end{prop}
\proof The proposition will follow from the following:

For each point-vertex $v$ in $\L$ there is a line $\la$ in $\R\P^2$ and for 
each line-vertex $v\in \L$ there is a line $\la'$ in $(\R\P^2)^{\vee}$ 
so that: 

$\phi(v)\notin \la$ for all $\phi\in  BR_0(\a, \R\P^2)$ (if $v$ 
is a point-vertex) and $\phi(v)\notin \la'$ for all 
$\phi\in  BR_0(\a, \R\P^2)$ (if $v$ is a line-vertex).  

To prove this property recall (see \cite{KM6}) that $\a$ is obtained 
from ``elementary''  arrangements for the addition and multiplication via 
fiber sums. Thus it is enough to verify the above property for the 
arrangements $C_A, C_M$ 
for the addition and multiplication that are described in \cite{KM6}. The 
verification is straightforward and is left to the reader.  \qed   

\medskip
Now we identify the moduli space of spherical linkages  
$\M(\L, {\mathbb S}^2)$ with a moduli space of Euclidean 
linkages in $\R^3$ as follows: 

\medskip 
Add an extra vertex $v_0$ to the graph $\L$ and connect it to each 
vertex of $\L$ by edge of the unit length. Modify the other 
side-lengths as follows:
$$
\ell'(e):= \sqrt{2 - 2\cos(\ell(e))}, \quad e\in \E(\L)
$$
Let $\L'$ be the resulting metric graph with the distinguished set of 
vertices $[\p(T)- \{v_{11}\}]\cup \{v_0\}$. Define the configuration space 
$$
C(\L', \R^3):= \{\psi: \V(\L') \to \R^3 : |\psi(v)- \psi(w)|^2= 
\ell'[vw]^2\}
$$
Again is is clear that 
$$
\M(\L', \R^3):= \{ \psi\in C(\L', \R^3): \psi(v_0)=(0,0,0), $$
$$
\hbox{and the same normalization on $\p(T)- \{v_{11}\}$ as we used for~~} 
\M(\L, {\mathbb S}^2)\} 
$$
is a real-algebraic set which is a cross-section for the action of 
$Isom(\R^3)$ on $C(\L', \R^3)$. Obviously we have an isomorphism 
$$
\M(\L, {\mathbb S}^2)\cong \M(\L', \R^3)
$$ 
of real-algebraic sets. We  let $\M_0(\L', \R^3)$ be the subset of 
$\M(\L', \R^3)$ corresponding to $\M_0(\L, \R\P^2)$ under the isomorphism
$$
\M(\L, \R\P^2)\cong \M(\L, {\mathbb S}^2)\cong \M(\L', \R^3)
$$

Thus, as a corollary of Theorem \ref{F} we obtain the following:

\begin{thm}
Let $S$ be a compact real algebraic set defined over $\Z$. Then there are 
abstract linkages $\L, \L'$ so that: 

\smallskip
(1) $\M_0(\L, \R\P^2)$ is entire rationally isomorphic to $S$. 

(2) $\M_0(\L', \R^3)$ is an (analytically) trivial entire rational 
covering of $S$. 

\smallskip
\noindent Both $\M_0(\L, \R\P^2)$, $\M_0(\L', \R^3)$ are Zariski open and 
closed subsets in the moduli spaces $\M(\L, \R\P^2)$, $\M(\L', \R^3)$ respectively. 
\end{thm}

%\begin{cor}
%Suppose that $M$ is a smooth compact manifold. Then there are linkages 
%$\L, \L'$ so that $M$ is diffeomorphic to unions of components in 
%$\M(\L, \R\P^2)$, $\M(\L', \R^3)$. 
%\end{cor}

\section{A brief history of ``Kempe's theorem''}
\label{hist}

This story began with the invention of the steam engine by 
Newcomen in 1722. 
One problem that appeared naturally was to transform a periodic linear 
motion (of the ``input'' vertex) to a circular motion  (of the ``output'' vertex).  The ``parallelogram'' invented 
by Watt in  the late 18-th century gave an approximate solution to 
this problem. The ``input'' motion was 
not exactly linear, however the input vertex traces a curve with a  
point of zero curvature, hence the output approximates a straight line 
up to the 2-nd order. After discovery in the first half of the 19-th century 
of several ``unsolvable'' geometric problems (like squaring a circle, etc.), 
 for a while it was a common opinion that the problem of transforming linear 
to circular motion also has no exact solution. This opinion was 
shared for instance by Chebyshev who after thinking about this problem  
introduced {\em Chebyshev polynomials}, partial motivation for which was 
the optimal approximate solution of the problem. 

This was the situation until 1864 when French navy officer Peaucellier 
published a letter \cite{Po1} where he claimed a positive solution, 
without giving any details\footnote{It seems that in 
1860-s Peaucellier explained his solution to some other people, cf. 
\cite{Mann}, so his letter \cite{Po1} was probably not a hoax. However 
\cite{Mann} contains only the title so we can not be sure if Mannheim 
really knew construction of the inversor.}.  

There are several opinions on what happened next (this caused a serious  
controversy between Russian and French-British mathematical schools in the late 
19-th century). In 1871 Lippman Lipkin\footnote{That time Lipkin was a 
graduate student of Chebyshev. Lipkin had died in 1875 at the age of 25 from 
the smallpox.} published the first detailed solution \cite{Lipkin}. 
Two years later (in 1873) Peaucellier published a paper \cite{Po2} 
which also contained a detailed solution  (the {\em Peaucellier inversor}) 
identical to Lipkin's.  Immediately after that several other ways to ``draw 
a straight line'' were discovered \cite{Hart}, \cite{Kempe2}. 
As far as applications are concerned it turned 
out that all the mechanisms that transform linear motion to circular  are 
too complicated to be used instead of Watt's parallelogram, invention of 
efficient lubricants had closed the problem. The only practical application 
of the inversor we are aware of was in air engines which ventilated 
the British parliament in 1870-1880-s (see \cite[Page 182]{W}). 

The rest of the story is mostly pure mathematics. In 1875 A.~B.~Kempe 
published \cite{Kempe} where (in the present terminology) he outlined a 
proof of the following theorem analogous to Theorem \ref{kem}: 

\begin{thm}
\label{Kempthm}
Suppose that $S\subset \R^2$ is an algebraic curve, $\O\in S$. Then 
there exists an abstract closed $\C$-functional linkage $\L$, a Zariski 
closed algebraic subset $C\subset \M(\L)$ (which is a union of 
irreducible components) and an closed \footnote{In the classical 
topology.} \nbd $U$ of $\O$ in $S$ so that the restriction of the 
input map $p$ to $C$ is onto $U$. 
\end{thm}

\begin{rem}
However, if one follows Kempe's arguments, 
$C$ is not open in $\M(\L)$, $U\ne S$ (even if $S$ is compact) 
and the mapping $p: C \to U$ is not a trivial covering.  
\end{rem}

Versions of Kempe's proof were reproduced in a number of places 
(see for instance \cite{B}), however (as far as we know) 
even the assertion was not made precise and details of the proof 
were not given. Recently several (written) attempts were made to improve 
Theorem \ref{Kempthm}, i.e. to make the subset $C$ open and $U= S$ 
(see \cite{HJW}) and the projection $p|C$ injective (see \cite{JS}),  
however, as far as we can tell, they were unsuccessful. 
Finally there was a work of W.~Thurston on this subject that 
we have discussed in the Introduction.

%\newpage

\noindent Michael Kapovich: Department of Mathematics, University of 
Utah,  Salt Lake City, UT 84112, USA ; kapovich$@$math.utah.edu

\smallskip
\noindent John J. Millson: Department of Mathematics, University of 
Maryland, College Park, MD 20742, USA ; jjm$@$math.umd.edu

\end{document}